\documentclass[12pt,twoside]{article}

\usepackage{graphicx}
\usepackage{amsfonts,amssymb,eucal,amsmath}
\usepackage{booktabs}

\newtheorem{CC}{Corollary}
\newtheorem{TT}{Theorem}
\newtheorem{LL}{Lemma}
\newtheorem{DD}{Definition}

\begin{document}

\setcounter{equation}{0}
 \thispagestyle{empty}

\begin{center}
{\Large\bf Vertex Structure \\ [4mm]
of Master Corner Polyhedra}\\[5mm]

{\large Vladimir A. Shlyk}\\%[3mm]

\end{center}

\renewcommand{\baselinestretch}{1}

\vspace*{5mm}
{\small \noindent
This paper focuses on vertices of the master corner polyhedra $P(G,g_0),$
the core of the group-theoretical approach to integer linear programming.
We introduce two combinatorial operations that transform each vertex of $P(G,g_0)$ to adjacent ones.
This implies that for any $P(G,g_0),$ there exists a subset of basic vertices, we call them support vertices,
from which all others can be built.
The class of support vertices is proved to be invariant under the automorphism group of $G,$
so this basis can be further reduced to a subset of pairwise non-equivalent support vertices.
Among other results, we characterize irreducible points of the master corner polyhedra,
establish relations between an integer point and the nontrivial facets that pass through it,
construct complete subgraphs of the graph of $P(G,g_0),$ and show that these polyhedra are of diameter~2.}

\section{Introduction}

The group-theoretical approach is one of the main approaches to the integer linear programming (ILP) \cite{NemWol99}.
It was originated by R. Gomory in the 1960s, cf. \cite{Gom65}, \cite{Gom69}. When applied to an ILP problem, this approach constructs its relaxation,
an optimization problem over a finite Abelian group $G.$ The polyhedron of its feasible solutions $P(G,H,g_0)$ is the convex hull
of solutions to the equation
$$
 \sum_{g\in H}t(g)g = g_0,
$$
where $H$ is some subset of $G,$ $g_0\in G.$
$P(G,H,g_0)$ is called the corner polyhedron or the Gomory polyhedron.
On certain conditions, the optimal solution to the initial ILP problem can be constructed from the optimal solution to the group optimization problem.
However, the main impact of the group-theoretical approach to the theory and practice of ILP is that
facets of corner polyhedra induce the most effective cuts for the initial optimization problems, including the mixed integer case.
These cuts are widely used within the branch-and-cut framework, some are implemented in commercial softwares, cf.~\cite{GloShe04-1page}.
As the "encyclopedia" of all corner polyhedra on $G$ with different $H$ and $g_0$ serves the master corner polyhedron
$P(G,g_0)=P(G,G^+,g_0),$ where $G^+=G\setminus\{\overline{0}\},$ $\overline{0}$ is the group zero element,
since its vertices and facets contain vertices and facets of all corner polyhedra $P(G,H,g_0)$ \cite{Gom69}.
So, $P(G,g_0)$ is the convex hull of the non-negative integer solutions $t=(t(g), g\in G^+)$ to the group equation
\begin{equation}\label{eq1.1}
 \sum_{g\in G^+}t(g)g = g_0
\end{equation}
and it lies in the $(D-1)$-dimensional space.

The computational potential of the group-theoretical approach remained unclear, if not questionable,
until the 1990s when a large variety of effective group cutting planes were elaborated.
This caused a revival of interest in the approach in general with the master corner polyhedron as its core object.
Since then, a great amount of research has been devoted to facets of $P(G,g_0)$ and generating effective cuts from them, see~\cite{Corn07},
\cite{GoJo03}, \cite{GoJoEvans03}, \cite{ArEvGoJo03}, \cite{RiDey10}.
By contrast, we are aware of only a few works, cf. \cite{Shh88sssr} for references, targeted on vertices during the 40 years period after the seminal Gomory's paper \cite{Gom69}.
To quote R. Gomory and E. Johnson \cite{GoJo04-1page},
"Understanding of these polyhedra, which one may well regard as the atoms of integer programming, is still at its beginning".

This paper focuses on the vertices of the master corner polyhedra though many its results can be transferred to the particular corner polyhedra.

We use three results due to R.~Gomory \cite{Gom69}.
The first, he proved that all vertices of $P(G,g_0)$ are irreducible points, where a point $t,$ a solution to (\ref{eq1.1}),
is irreducible if for any points $r=(r(g), g\in G^+)$ and $s=(s(g), g\in G^+),$ the conditions
\begin{equation}\label{eq.Irreduc}
 0\leq r(g), s(g)\leq t(g), \quad \sum_{g\in G^+}r(g)g = \sum_{g\in G^+}s(g)g,
\end{equation}
imply $r=s.$

The second, each automorphism $\varphi$ of the group $G$ transforms any vertex of $P(G,g_0)$ to a vertex of $P(G,\varphi(g_0)).$

The third is the subadditive description of facets of $P(G,g_0).$ The facets are canonically considered as inequalities
\begin{equation}\label{eq1.3}
 \sum_{g\in G^+}\pi(g)t(g) \geq \pi_0
\end{equation}
and are denoted by $(\pi,\pi_0),$ with $\pi=(\pi(g), g\in G^+)$ as their coefficient vectors.
Gomory proved that facets of $P(G,g_0)$ are of two types.
The trivial facets are inequalities $t(g)\geq 0,$ $g\in G^+,$ determining coordinate hyperplanes in $\mathbb{R}^{D-1}.$
The nontrivial facets are the inequalities $(\pi,\pi_0),$ $\pi_0>0,$ with $\pi$ a basic feasible solution to the system
\begin{subequations}\label{eq4:ab}
\begin{align}\label{eq4:a}
&\pi(g_0)=\pi_0,
\\ \nonumber
&\pi(g)+\pi(g_0-g)=\pi_0, &g\in G^+, g\neq g_0,
\\ \label{eq4:b}
&\pi(g_1)+\pi(g_2)\geq\pi(g_1+g_2), &g_1,g_2\in G^+,
\\ \nonumber
&\pi(g)\geq 0, &g\in G^+,
\end{align}
\end{subequations}
where (\ref{eq4:a}) is omitted if $g_0=\overline{0}.$

This paper is organized as follows. In the next section, we prove a geometric characterization of the irreducible points of
master corner polyhedra: these are such integer points of a $P(G,g_0)$ that cannot be expressed as a convex combination
of any two its integer points.
In Section 3, we show that coefficients of the nontrivial facets passing through a given integer vertex of $P(G,g_0)$
turn certain subadditive inequalities (\ref{eq4:b}) to equalities and construct some other points in these facets.
Section~4 is central to the paper.
We introduce two combinatorial $\mu$-operations and prove that, when applied to vertices of a $P(G,g_0),$ they result in adjacent vertices.
This implies that every master corner polyhedron can be determined by a subset of its basic vertices,
those from which all others can be built with the use of $\mu$-operations; we call them support vertices.
We prove a support vertex analog of Gomory's theorem on automorphisms: each automorphism $\varphi$ of $G$
transforms any support vertex of $P(G,g_0)$ to a support vertex of $P(G,\varphi(g_0)).$
Together with the concordance of the $\mu$-operations with automorphisms, this result provides a more detailed view of the
vertex structure of master corner polyhedra and leads to a description of their minimal (with respect to the tools at hand) vertex bases.
In Section~5, applying $\mu$-operations recursively, we construct some vertex sequences that generate complete subgraphs
of the graph of $P(G,g_0).$ As a consequence it follows that master corner polyhedra are of diameter 2.
Section 6 presents the conclusions of the work and discusses some directions for future study.

We denote by $V(G,g_0)$ the set of vertices of $P(G,g_0)$ and for $t\in P(G,g_0),$ set $G_t=\{g\in G^+|t(g)>0\}.$
The standard notations $Aut(G)$ and $Aut_h(G)$ are used respectively for the automorphism group of $G$
and the stabilizer subgroup of $Aut(G)$ at $h\in G:$
$$
Aut_h(G)=\{\varphi\in Aut(G)|\varphi(h)=h\}.
$$
We refer the reader to \cite{Her75} and \cite{Asch00} for the details on the group automorphisms and the group actions.

\section{Geometric characterization of irreducibility}

In this section, we prove a geometric characterization of irreducible points of master corner polyhedra,
thereby extending Theorem 2 in \cite{Gom69}. But first, we present several properties of irreducible points.

\begin{LL}\label{lem1-Irred}
Let $t\in P(G,g_0)$ be an irreducible point and $u=(u(g),g\in G^+),$ $u\neq t,$ have the components $0\leq u(g)\leq t(g),$ $g\in G^+.$ Then
\begin{itemize}
\item[{\rm(i)}] $\sum_{g\in G^+}u(g)g\notin G_t,$
\item[{\rm(ii)}] $\sum_{g\in G^+}u(g)g\neq \overline{0},$
\item[{\rm(iii)}] $\sum_{g\in G^+}u(g)g\neq g_0.$
\end{itemize}
\end{LL}

\noindent{\bf Proof.}
(i) is a straightforward consequence of irreducibility. To prove (ii), set $r=t$ and $s=(t(g)-u(g), g\in G^+).$
Then $r$ and $s$ satisfy (\ref{eq.Irreduc}),
which together with $r\neq s$ contradicts irreducibility of $t.$
To prove (iii), notice that otherwise the equality $\sum_{g\in G^+}(t(g)-u(g))g\!=\!\overline{0}$ would contradict (ii).
$\hfill\square$

\begin{TT}\label{Th1-Irred}
An integer point of $P(G,g_0)$ is irreducible if and only if it cannot be expressed as a convex combination
of any two its integer points.
\end{TT}

\noindent{\bf Proof.}
Gomory proved irreducibility of vertices by showing that any reducible point $t\in P(G,g_0)$
is a half-sum of two solutions to equation (\ref{eq1.1}), which implies that $t$ is a convex combination of two integer points of $P(G,g_0).$

So, it remains only to show that if an integer point $t\in P(G,g_0)$ is a convex combination of some integer points
$u,v\in P(G,g_0)$ then $t$ is reducible.
One can easily see that if $t,$ a point of the $n$-dimensional integer grid, is a convex combination of two grid points
$u$ and $v$ then $t$ is the half-sum of the two nearest to $t$ grid points $u'$ and $v'$ in the line segment $[u,v].$
The inclusions $u',v'\in P(G,g_0)$ hold by convexity. So, we can deal with $t=\frac12(u+v),$ $u,v\in P(G,g_0).$
Then
\begin{equation}\label{eq2.3}
 t(g)-u(g)=v(g)-t(g) ~~\mathrm{for ~all}~~ g\in G^+
\end{equation}
and we can construct integer points $r$ and $s$ as follows,
\begin{equation}\label{Th1:ProIrred}
 \begin{array}{lll}
&r(g)=t(g)-u(g), \quad &g\in H_u=\{g\in G^+\mid u(g)<t(g)\},
\\
&r(g)=0, \quad &g\in G^+\setminus H_v,
\\
&s(g)=t(g)-v(g), \quad &g\in H_v=\{g\in G^+\mid v(g)<t(g)\},
\\
&s(g)=0, \quad &g\in G^+\setminus H_v.
 \end{array}
\end{equation}
Notice that $H_u\cap H_v=\emptyset$ and if $H_u=\emptyset$ or $H_v=\emptyset$ is empty then $r=t$ or $s=t.$
By (\ref{eq2.3}) and (\ref{Th1:ProIrred}), $r$ and $s$ satisfy~(\ref{eq.Irreduc}).
However, $r\neq s,$ which means that $t$ is reducible and ends the proof. $\hfill\square$

Theorem \ref{Th1-Irred} clarifies why some irreducible points of master corner polyhedra are not vertices:
they are convex combinations of some more than two points. The following example presents such a point.

\vspace{2mm} \noindent
{\bf Example 1.}
Consider the master corner polyhedron $P(G_{16},\overline{15}),$ where $G_{16}$ is the cyclic group of order 16 with the elements $\overline{0},\overline{1},\overline{2},\ldots,\overline{15}.$ We will use $0^k$ for the sequence of $k$ zeroes.
The point
$$
t=(0,0,2,1,1,0^{10})
$$
belongs to $P(G_{16},\overline{15})$ as it satisfies the equality
$2\cdot\overline{3}+1\cdot\overline{4}+1\cdot\overline{5}=\overline{15}.$
One can check that it is irreducible by comparing pairwise the sums
$r(\overline{3})\cdot\overline{3}+r(\overline{4})\cdot\overline{4}+r(\overline{5})\cdot\overline{5}$
of $3\cdot2\cdot2-1=11$ integer points
$r=(r(g), g\in G_{16}^+),$ $0\leq r(g)\leq t(g),$ $r\neq(0^{15}).$
Since these sums are distinct, Theorem \ref{Th1-Irred} implies that $t$ is not a convex combination of any two integer points of $P(G_{16},\overline{15}).$ However $t$ is the convex combination of three points,
$$
t=\frac13(0,0,0,0,3,0^{10})+\frac13(0,0,1,3,0,0^{10})+\frac13(0,0,5,0,0,0^{10}),
$$
and thus is not a vertex.

\section{Integer points and nontrivial facets}

In this section, we obtain some relations between integer points of master corner polyhedra and their nontrivial facets.

\begin{TT}\label{Th2-Vert-Fac}
Let a vertex $t$ of the polyhedron $P(G,g_0)$ belong to its nontrivial facet $(\pi,\pi_0),$
and let $h=\sum_{g\in G^+}u(g)g$ with integer $u(g),$ $0\leq u(g)\leq t(g),$ $g\in G^+.$ Then
\begin{itemize}
\item[{\rm(i)}] the point $w=(w(g), g\in G^+)$ with the components $w(g)=t(g)-u(g)$ for $g\in G^+,$ $g\neq h,$ and $w(h)=t(h)+1$
belongs to all nontrivial facets of $P(G,g_0)$ that contain$~t;$
\item[{\rm(ii)}] coefficients of the facet $(\pi,\pi_0)$ satisfy the relation
$$
\pi(h)=\sum_{g\in G_t}u(g)\pi(g).
$$
\end{itemize}
\end{TT}

\noindent{\bf Proof.} To prove (i), observe that $w$ is a solution to the group equation (\ref{eq1.1}), hence $w\in P(G,g_0).$
As $(\pi,\pi_0)$ is a facet, $w$ satisfies the inequality
$$
\sum_{g\in G^+}w(g)\pi(g)\geq \pi_0.
$$
It follows from the subadditivity condition (\ref{eq4:b}) that 
$$
\sum_{g\in G^+}u(g)\pi(g)=\sum_{g\in G_t}u(g)\pi(g)\geq \pi(h),
$$ 
while Lemma~\ref{lem1-Irred} (i) implies $t(h)=0.$ Therefore,
\begin{equation*}
\begin{aligned}
\sum_{g\in G^+}w(g)\pi(g)
&
 =\!\sum\limits_{\begin{array}{c}
 \scriptstyle~\\[-7.5mm]
 \scriptstyle g\in G^+ \\[-1mm]
 \scriptstyle g\neq h \\[-1mm]
 \end{array}} \!(t(g)-u(g))\pi(g)+(t(h)+1)\pi(h)
\\
&
\leq\!\sum\limits_{\begin{array}{c}
 \scriptstyle~\\[-7.5mm]
 \scriptstyle g\in G^+ \\[-1mm]
 \scriptstyle g\neq h \\[-1mm]
 \end{array}} \!t(g)\pi(g)-
 \!\sum\limits_{\begin{array}{c}
 \scriptstyle~\\[-7.5mm]
 \scriptstyle g\in G^+ \\[-1mm]
 \scriptstyle g\neq h \\[-1mm]
 \end{array}} \!u(g)\pi(g)
 +
 \!\sum\limits_{\begin{array}{c}
 \scriptstyle~\\[-7.5mm]
 \scriptstyle g\in G^+ \\[-1mm]
 \scriptstyle g\neq h \\[-1mm]
 \end{array}} \!u(g)\pi(g)
\\
&
 =\sum_{g\in G^+}t(g)\pi(g)=\pi_0.
\end{aligned}
\end{equation*}
The two opposite inequalities imply the equality
\begin{equation}\label{eq3.1}
\sum_{g\in G^+}w(g)\pi(g)=\pi_0,
\end{equation}
which proves (i). One can notice that~(\ref{eq3.1}) holds in the only case of
$\pi(h)=\sum_{g\in G_t}u(g)\pi(g).$ This proves (ii) and completes the proof of the theorem. $\hfill\square$

\section{Support vertices}

This section concerns several topics related to some special vertices of the master corner polyhedron.
We introduce two combinatorial operations that can be applied to the most of its integer points and prove that they transform vertices to adjacent vertices. 
This leads to existence of a subset of vertices of each master corner polyhedron,
called support vertices, that do not result from any other vertex by these operations.
Then we study the structure of the orbit partition of the set of vertices of the $P(G,g_0)$ under the action of the automorphism group of $G.$
We prove that some orbits fully consist of support vertices and describe bases of the set $V(G,g_0).$

\subsection{$\mu$-operations and support vertices}

Let us define two combinatorial operations applicable to some integer points of the $P(G,g_0).$
We call them $\mu$-operations.

\textbf{Operation} $\mu_{h,f}.$ Let $t$ be an integer point of $P(G,g_0)$ and let $h,f\in G_t;$ for certainty, let $t(h)\leq f(h).$
Build the point $s=\mu_{h,f}(t)$ with the components
\begin{equation*}
\begin{aligned}
&
 s(h)=0, \\
&
s(f)=t(f)-t(h), \\
&
s(h+f)=t(h+f)+t(h), \\
&
s(g)=t(g), \quad \quad \quad g\in G^+, ~g\neq h,f,h+f.
\end{aligned}
\end{equation*}

\textbf{Operation} $\mu_{h}.$ Let $t$ be an integer point of $P(G,g_0)$ and let $h\in G_t$ satisfy $t(h)>1.$
Build the point $s=\mu_h(t)$ with the components
\begin{equation*}
\begin{aligned}
&
 s(h)=0, \\
&
s(t(h)h)=t(t(h)h)+1, \\
&
s(g)=t(g), \quad \quad \quad g\in G^+, ~g\neq h,t(h)h.
\end{aligned}
\end{equation*}

\begin{TT}\label{Th-MuAdjVert}
Let $t$ be a vertex of $P(G,g_0)$ and the operation $\mu_{h,f}$ with some $h,f\in G_t$
(respectively, $\mu_{h}$ with some $h\in G_t$) be applicable to $t.$
Then $\mu_{h,f}(t)$ (respectively, $\mu_{h}(t)$) is a vertex of $P(G,g_0)$ adjacent to $t.$
\end{TT}

\noindent{\bf Proof.}
We will prove the theorem for the case of $\mu_{h,f}$ since the case of $\mu_{h}$ can be considered similarly.
At first, prove that the point $s=\mu_{h,f}(t)$ is a vertex of $P(G,g_0)$ provided $t$ is a vertex.
It is an easy exercise to check that $s\in P(G,g_0).$ Assume $s$ is not a vertex.
Then $s$ is a convex combination of some $k\geq 2$ integer points $s_j,$ $j=1,2,\ldots,k,$ that solve equation (\ref{eq1.1}):
$s=\sum_{J=1}^k\lambda_js_j,$ $\sum_{J=1}^k\lambda_j=1,$ $\lambda_j>0.$
It follows from $s(h)=0$ that $s_j(h)=0$ for all $j.$
Define integer points $t_j,$ $j=1,2,\ldots,k,$ by setting
 $$
 \begin{array}{lll}
&t_j(h)=s_j(h+f); \quad &t_j(f)=s_j(h+f)+s_j(f);
\\
&t_j(h+f)=0; \quad &t_j(g)=s_j(g),~g\in G^+\setminus \{h,f,h+f\}
 \end{array}
 $$
and check that all $t_j\in P(G,g_0):$
\begin{equation*}
\begin{aligned}
\sum_{g\in G^+}t_j(g)g&=t_j(h)h+t_j(f)f+t_j(h+f)(h+f)+\!\sum_{g\in G^+\setminus \{h,f,h+f\}}\!t_j(g)g
\\
&
=s_j(h+f)h+s_j(h+f)f+s_j(f)f+\!\sum_{g\in G^+\setminus \{h,f,h+f\}}\!s_j(g)g
\\
&
=s_j(h+f)(h+f)+s_j(f)f+\!\sum_{g\in G^+\setminus \{h,f,h+f\}}\!s_j(g)g
\\
&
=\sum_{g\in G^+}s_j(g)g=g_0.
\end{aligned}
\end{equation*}

By Lemma \ref{lem1-Irred} (i), $t(h\!+\!f)\!=\!0.$ Using this equality, show that $\sum_{j=1}^k\!\lambda_jt_j\!=\!t:$
\begin{equation*}
\begin{aligned}
&
\sum_{j=1}^k\lambda_jt_j(h)=\sum_{j=1}^k\lambda_js_j(h+f)=s(h+f)=t(h+f)+t(h)=t(h),
\\
&
\sum_{j=1}^k\lambda_jt_j(f)=\sum_{j=1}^k\lambda_js_j(h+f)+\sum_{j=1}^k\lambda_js_j(f)=s(h+f)+s(f)
\\
&~~~~~~~~~~~~~~~
=t(h+f)+t(h)+t(f)-t(h)=t(f),
\\
&
\sum_{j=1}^k\lambda_jt_j(h+f)=0=t(h+f),
\\
&
\sum_{j=1}^k\lambda_jt_j(g)=t(g), \quad~~ g\in G^+\setminus \{h,f,h+f\}.
\end{aligned}
\end{equation*}

So, we have obtained that $t$ admits a convex representation via $k$ points of $P(G,g_0),$
however this contradicts $t$ being a vertex. Therefore, $s$ is a vertex of $P(G,g_0).$

Now prove that $s$ is adjacent to $t.$
As every vertex $v$ of a full-dimensional polyhedron $P\subset\mathbb{R}^n$ can be defined by a collection $\mathcal{F}_v$ of
$n$ linearly independent facets of $P$ passing through $v,$ it is sufficient to find some collections $\mathcal{F}_t$ and $\mathcal{F}_s$
of $D-1$ linearly independent facets of $P(G,g_0)$ each, passing respectively through $t$ and $s$
and differing by only one facet. Then their common facets define the edge $(t,s)$ of $P(G,g_0).$

First, include into $\mathcal{F}_t$ all trivial facets $t(g)\ge0,$ $g\in G^+\setminus G_t,$ of $P(G,g_0),$ which obviously contain $t.$
By Lemma \ref{lem1-Irred} (i), $t(h+f)\ge0$ is one of these facets. Next, add into $\mathcal{F}_t$ the necessary amount
($=D-1-|G^+\setminus G_t|=|G_t|$) of nontrivial facets $(\pi,\pi_0)$ passing through $t$
and such that all facets in $\mathcal{F}_t$ be linearly independent. As $t$ is a vertex such facets exist.
The coefficient matrix $M_t$ of the facets in $\mathcal{F}_t$ contains the identity submatrix disposed
in the rows corresponding to the trivial facets and the columns indexed by $j\in G^+\setminus G_t;$
notice that the $(h+f)$-th column is one of these. Therefore, the nontrivial facets in $\mathcal{F}_t$
are linearly independent on the columns $j\in G_t.$

Now, build the collection $\mathcal{F}_s.$ Include into it all facets from $\mathcal{F}_t$ except $t(h+f)\!\ge\!0,$
instead of which use the trivial facet $t(h)\!\ge\!0.$ The trivial facets in $\mathcal{F}_s$ contain $s$ by its construction,
while the nontrivial facets contain $s$ by Theorem~\ref{Th2-Vert-Fac} (i).
So, the coefficient matrix $M_s$ of the facets in $\mathcal{F}_s$ contains the identity submatrix in the columns
$j\in J_1=((G^+\setminus G_t)\setminus\{h+f\})\cup\{h\}.$
The nontrivial facets in $\mathcal{F}_t$ and $\mathcal{F}_s$ are the same and, by Theorem \ref{Th2-Vert-Fac} (ii), their coefficients
satisfy the equality $\pi(h)+\pi(f)=\pi(h+f).$ This yields linear independence of the nontrivial facets in $\mathcal{F}_s$ on the columns
indexed by $j\in J_2=(G_t\setminus\{h\})\cup\{h+f\}.$ As $J_1\cap J_2=\emptyset$ the facets in $\mathcal{F}_s$ are linearly independent,
so $\mathcal{F}_s$ is the collection we strived to obtain. The theorem is proved. $\hfill\square$

\begin{DD}
We call a vertex $t$ of a master corner polyhedron $P(G,g_0)$ a support vertex if $t$ does not result
from any other vertex of this polyhedron with the use of any $\mu$-operation.
\end{DD}

The inequality $\sum_{g\in G^+}s(g)g < \sum_{g\in G^+}t(g)g,$ provided $s=\mu_{h,f}(t)$ or $s=\mu_h(t),$ 
implies existence of support vertices of any $P(G,g_0).$ They are of special importance because by Theorem~\ref{Th-MuAdjVert},
they form a basis for $V(G,g_0),$ since every other vertex can be build from some support vertex by recursive application of some $\mu$-operations.

Denote the set of support vertices of $P(G,g_0)$ by $S(G,g_0).$
The next example presents support vertices of $P(G_6,\overline{3}).$

\vspace{2mm} \noindent
{\bf Example 2.}
Let us continue using notation from Example~1. We find from Table~1 in~\cite{Gom69} that $P(G_6,\overline{3})$ has 7 vertices:
\begin{equation*}
\begin{aligned}
& t_1=(3,0,0,0,0), ~~~~~~~~~~
& t_2=(1,1,0,0,0),
\\
& t_3=(0,0,1,0,0), ~~~~~~~~~~
& t_4=(1,0,0,2,0),
\\
& t_5=(0,2,0,0,1), ~~~~~~~~~~
& t_6=(0,0,0,1,1),
\\
& t_7=(0,0,0,0,3)
\end{aligned}
\end{equation*}
and four nontrivial facets. The support vertices are only $t_1,$ $t_4,$ $t_5,$ $t_7$ since
$t_2=\mu_4(t_4),$ $t_3=\mu_{1,2}(t_2),$ and $t_6=\mu_{1,4}(t_4),$ though there are more ways to obtain $t_2,$ $t_3,$ $t_6.$
One can observe all vertices of this polyhedron together with the $\mu$-operations acting on them in Figure~1 below.

The vertex $t_4$ belongs to the nontrivial facets
$(\pi,\pi_0)=$
\begin{equation*}
\begin{aligned}
& ((1,0,1,0,1,1),1), ~~~~~~~~~~
& ((1,2,3,1,2,3),3)
\end{aligned}
\end{equation*}
and the trivial facets $t(2)\ge0,$ $t(3)\ge0,$ $t(5)\ge0.$ The vertex $t_2$ belongs to the same nontrivial facets, as well as to the facets
\begin{equation*}
\begin{aligned}
& ((2,1,3,2,1,3),3), ~~~~~~~~~~
& ((1,2,3,2,1,3),3).
\end{aligned}
\end{equation*}
So, in fact, $t_2$ belongs to all nontrivial facets of $P(G_6,\overline{3})$ and to the trivial facets
$t(3)\ge0,$ $t(5)\ge0,$ and $t(4)\ge0$ instead of $t(2)\ge0.$

The vertex $t_3=\mu_{1,2}(t_2)$ belongs to all four nontrivial facets and the trivial facets $t(4)\ge0,$ $t(5)\ge0,$ and $t(1)\ge0,$ $t(2)\ge0.$
Each of the last two can be considered as substituting $t(3)\ge0.$

%\newpage
\subsection{Automorphisms and support vertices}

R. Gomory proved that each automorphism $\varphi$ of the group $G$ transforms any vertex $t$ of $P(G,g_0)$ to a vertex
\begin{equation}\label{eq8}
\overline{t}=\{\overline{t}(g),g\in G^+\}=\{t(\varphi^{-1}(g)),g\in G^+\}
\end{equation}
of $P(G,\varphi(g_0)),$ see \cite{Gom69}, the Corollary following Theorem 14.

Let $\mathcal{V}(G)$ denote the set of vertices of all master corner polyhedra on $G,$ $\mathcal{V}(G)=\cup_{g_0\in G}V(G,g_0).$
Then (\ref{eq8}) defines the binary function
\begin{equation}\label{eq9}
\mathcal{V}(G)\times Aut(G)\rightarrow \mathcal{V}(G)~:~(t,\varphi)\mapsto t\cdot\varphi=\overline{t}.
\end{equation}

Since $(\varphi\sigma)^{-1}=\sigma^{-1}\varphi^{-1},$ this function satisfies two conditions
\begin{equation*}
\begin{aligned}
&t\cdot(\varphi\sigma)=(t\cdot\varphi)\cdot\sigma ~~\mathrm{for~all}~~ t\in \mathcal{V}(G) ~~\mathrm{and}~~ \varphi,\sigma\in Aut(G),
\\
&t\cdot\varepsilon=t ~~~~~~~~~~~~~~~~~\mathrm{for~all}~~ t\in \mathcal{V}(G) ~~\mathrm{and}~~ \varepsilon, ~\mathrm{the~identity~automorphism~of}~~ G.
\end{aligned}
\end{equation*}

This means that~(\ref{eq9}) defines a (right) group action of $Aut(G)$ on $\mathcal{V}(G):$
each automorphism $\varphi\in Aut(G)$ is represented as the vertex transformation sending any $t\in \mathcal{V}(G)$
to $\overline{t}\in \mathcal{V}(G).$ By (\ref{eq8}), the vertex $\overline{t}=t\cdot\varphi$ has the same as $t$ though permuted components.

As Gomory's theorem specifies the polyhedron $P(G,h),$ $h\in G,$ of which $t\cdot\varphi$ is a vertex,
it particularly asserts that $V(G,g_0)$ is invariant under the action of $Aut_{g_0}(G)$ on $V.$
We will call vertices $t$ and $s$ of some $P(G,g_0)$ equivalent if $s=t\cdot\varphi$ for some $\varphi\in Aut_{g_0}(G).$
Note that then by (\ref{eq8}), $t=s\cdot\varphi^{-1},$ $\varphi^{-1}\in Aut_{g_0}(G).$

The next lemma states that up to a slight change in the the group elements $h,f$ the $\mu$-operations commute with automorphisms of $G$
acting on $\mathcal{V}(G).$

\begin{LL}\label{Lem_Mu-Phi}
Let $t$ be a vertex of $P(G,g_0)$ such that some operation $\mu_{h,f}$ with $h,f\in G_t,$ or some operation $\mu_h$ with $h\in G_t,$
be applicable to $t,$ and let $\varphi\in Aut(G).$ Then
\begin{equation}\label{commutHF}
(\mu_{h,f}(t))\cdot\varphi=\mu_{\varphi(h),\varphi(f)}(t\cdot\varphi)
\end{equation}
or, respectively,
\begin{equation}\label{eqCommutH}
(\mu_h(t))\cdot\varphi=\mu_{\varphi(h)}(t\cdot\varphi).
\end{equation}
\end{LL}

\noindent{\bf Proof.}
We will prove only (\ref{commutHF}) since (\ref{eqCommutH}) can be proved similarly. Denote $s=\mu_{h,f}(t)$ and $\overline{t}=t\cdot\varphi.$
By Theorem \ref{Th-MuAdjVert} and Gomory's theorem, $s$ and $\overline{t}$ are some vertices of $P(G,g_0)$ and $P(G,\varphi(g_0)),$ respectively,
and we must prove that
\begin{equation}\label{phi(s)}
s\cdot\varphi=\mu_{\varphi(h),\varphi(f)}(\overline{t}).
\end{equation}
We show this equality component-wise using the rule (\ref{eq8}) and the definition of $\mu_{h,f}.$

\begin{equation*}
\begin{aligned}
(s\cdot\varphi)(\varphi(h))=&s(\varphi^{-1}\varphi(h))=s(h)=0=\big{[}\mu_{\varphi(h),\varphi(f)}(\overline{t})\big{]}(\varphi(h)),
\\
(s\cdot\varphi)(\varphi(f))=&s(\varphi^{-1}\varphi(f))=s(f)=t(f)-t(h)
\\
=&\overline{t}(\varphi(f))-\overline{t}(\varphi(h))=\big{[}\mu_{\varphi(h),\varphi(f)}(\overline{t})\big{]}(\varphi(f)),
\\
(s\cdot\varphi)(\varphi(h+f))=&s(\varphi^{-1}\varphi(h+f))=s(h+f)=t(h)=\overline{t}(\varphi(h))=
\\
=&\big{[}\mu_{\varphi(h),\varphi(f)}(\overline{t})\big{]}(\varphi(h))=\big{[}\mu_{\varphi(h),\varphi(f)}(\overline{t})\big{]}(\varphi(h+f),
\\
(s\cdot\varphi)(\varphi(g))=&s(\varphi^{-1}\varphi(g))=s(g)=t(g)=\overline{t}(\varphi(g))=\big{[}\mu_{\varphi(h),\varphi(f)}(\overline{t})\big{]}(\varphi(g))
\\
&\text{for} ~g\in G^+, g\neq h,f,h+f.
\end{aligned}
\end{equation*}
Thus (\ref{phi(s)}), (\ref{commutHF}), and the lemma are proved. $\hfill\square$

Paraphrasing the lemma brings us to two assertions advantageous to a better view of the vertex set structure.
Now, we will omit the group elements $h$ and $f$ determining the $\mu$-operation in use, though leave the index $\varphi$ to determine the new
$\mu$-operation $\mu_{\varphi(h),\varphi(f)}$ or $\mu_{\varphi(h)}$ induced by an automorphism $\varphi.$

\begin{CC}\label{Cor_Mu-Phi}
\begin{itemize}
\item[{\rm(i)}]
Let $t$ and $s$ be some vertices of $P(G,g_0)$ and let $\varphi\in Aut(G).$ Then $s=\mu(t)$ for some $\mu$-operation $\mu$ if and only if
$s\cdot\varphi=\mu_{\varphi}(t\cdot\varphi).$
\item[{\rm(ii)}]
Let $t$ and $s$ be some vertices of respectively $P(G,g_0)$ and $P(G,\varphi(g_0)),$ where $\varphi\in Aut(G),$
and let some $\mu$-operation $\mu$ be applicable to $t.$
Then $s=t\cdot\varphi$ if and only if $\mu_{\varphi}(s)=\mu(t)\cdot\varphi.$
\end{itemize}
\end{CC}

Note that (ii) does not preclude the equality $\mu_{\varphi}(s)=\mu(t)$ for $s\neq t.$
One can observe this case in Figure~1 for the vertices $t_1=(3,0,0,0,0),$ $t_7=t_1\cdot\varphi_5=(0,0,0,0,3),$
and $t_3=(0,0,1,0,0)=\mu_{\overline{1}}(t_1)=\mu_{\overline{5}}(t_7).$

The next theorem provides the support vertex analog of Gomory's theorem.
\begin{TT}\label{Th_SupV-Auto}
For any support vertex $t$ of some $P(G,g_0)$ and any automorphism $\varphi$ of $G,$ $t\cdot\varphi$ is a support vertex of $P(G,\varphi(g_0)).$
\end{TT}

\noindent{\bf Proof.}
Assume $\overline{t}=t\cdot\varphi\in V(G,\varphi(g_0))$ is not support, then $\overline{t}=\mu(\overline{r})$
for some $\mu$ and $\overline{r}\in V(G,\varphi(g_0)).$ By Corollary \ref{Cor_Mu-Phi} (i),
$t=\overline{t}\cdot\varphi^{-1}=\mu_{\varphi^{-1}}(\overline{r}\cdot\varphi^{-1}),$
where $\overline{r}\cdot\varphi^{-1}\in V(G,g_0).$ This contradicts $t\in S(G,g_0),$ so the theorem is proved. $\hfill\square$

\begin{CC}\label{Cor_Sup-Vert}
\begin{itemize}
\item[{\rm(i)}]
The set of support vertices of any $P(G,\overline{0})$ is invariant under $Aut(G)$ acting on $\mathcal{V}(G).$
\item[{\rm(ii)}]
The set of support vertices of any $P(G,g_0)$ with $g_0\neq \overline{0}$ is invariant under $Aut_{g_0}(G)$ acting on $V(G,g_0).$
\end{itemize}
\end{CC}

\noindent{\bf Proof.} The corollary follows from the theorem since for any $\varphi\in Aut(G)$ and $g\in G,$
$\varphi(g)=\overline{0}$ if and only if $g=\overline{0}.$ $\hfill\square$

\vspace{2mm} \noindent
{\bf Example 3.}
Now, we can consider support vertices
\begin{equation*}
\begin{aligned}
&t_1=(3,0,0,0,0), ~~~~~~~&t_4=(1,0,0,2,0),
\\
&t_5=(0,2,0,0,1), ~~~~~~~
&t_7=(0,0,0,0,3).
\end{aligned}
\end{equation*}
of $P(G_6,\overline{3})$ found in Example 2 more thoroughly.
There is only one non-identity automorphism $\varphi_5$ of the cyclic group $G_6,$
which maps each group element $\overline{r}$ to $\overline{5r\!\!\pmod{6}},$ $r=0,1,\ldots,5.$
Luckily, $\varphi_5$ leaves the right-hand-side element $\overline{3}$ fixed, so $\varphi_5\in Aut_{\overline{3}}(G_6).$
Thus, by Corollary 3 (ii), $\varphi_5$ transforms each support vertex of $P(G_6,\overline{3})$ to a support vertex of this polyhedron
(which may be the same vertex). One can see that
$$
t_1\cdot\varphi_5=t_7, ~~~~~~~~~t_4\cdot\varphi_5=t_5,
$$
and vice versa. As a result, the set $\{t_1,t_4\}$ of two support vertices can be regarded as a basis of $V(G_6,\overline{3}),$
as well as three more pairs of support vertices, see Figure~1.
If we take into account only automorphisms of $G_6$ but not the $\mu$-operations,
we will come to a basis consisting of four vertices, $\{t_1,t_2,t_3,t_4\}$ as an instance.

\begin{center}
\linethickness{0.15mm}
\begin{picture}(380,150) \label{fig.1} %\caption{Vertices of $P(G_6,\overline{3})$ with $\mu$-operations and automorphisms}\label{figure1}
\scriptsize
\put(0,0){\line(1,0){380}} \put(0,140){\line(1,0){380}} \put(0,0){\line(0,1){140}} \put(380,0){\line(0,1){140}}
\put(70,70){\begin{picture}(40,40)
% Orbity
\put(0,0){\oval(10,100)} \put(-5,-5){\vector(0,1){5}} \put(5,5){\vector(0,-1){5}}
        \put(-16,-10){$\varphi_5$} \put(8,10){$\varphi_5$} % Orbita 1
\put(80,0){\oval(10,100)} \put(75,-5){\vector(0,1){5}} \put(85,5){\vector(0,-1){5}} % Orbita 2
        \put(64,-10){$\varphi_5$} \put(89,10){$\varphi_5$}
\qbezier(163,5)(175,30)(160,30) \qbezier(157,5)(145,30)(160,30) \put(157,30){\vector(1,0){6}} % Orbit 3 -- Petlya
        \put(155,35){$\varphi_5$}
\put(240,0){\oval(10,100)} \put(235,-10){\vector(0,1){5}} \put(245,5){\vector(0,-1){5}} % Orbita 4
        \put(224,-18){$\varphi_5$} \put(248,10){$\varphi_5$}

% Tochki na orbitah
\put(0,50){\circle*{6}} \put(0,-50){\circle*{6}}       \put(0,50){\circle{8}} \put(0,-50){\circle{8}}
\put(80,50){\circle*{6}} \put(80,-50){\circle*{6}}
\put(160,0){\circle*{6}}
\put(240,50){\circle*{6}} \put(240,-50){\circle*{6}}   \put(240,50){\circle{8}} \put(240,-50){\circle{8}}
% Mu-strelki
\put(8,50){\vector(1,0){64}}     \put(8,-50){\vector(1,0){64}}     \put(40,55){$\mu_4$}       \put(40,-58){$\mu_2$}
\put(8,45){\vector(3,-4){64}}    \put(8,-45){\vector(3,4){64}}     \put(23,30){$\mu_{1,4}$}   \put(23,-33){$\mu_{5,2}$}
\put(88,46){\vector(3,-2){64}}   \put(88,-45){\vector(3,2){64}}    \put(117,30){$\mu_{1,2}$}  \put(117,-33){$\mu_{5,4}$}
\put(232,46){\vector(-3,-2){64}} \put(232,-45){\vector(-3,2){64}}  \put(196,30){$\mu_1$}      \put(197,-33){$\mu_5$}
% Podpisi
\put(-55,58){$t_4=(1,0,0,2,0)$}       \put(-55,-64){$t_5=(0,2,0,0,1)$}
\put(75,58){$t_2=(1,1,0,0,0)$}        \put(75,-64){$t_6=(0,0,0,1,1)$}
\put(235,58){$t_1=(3,0,0,0,0)$}       \put(235,-64){$t_7=(0,0,0,0,3)$}
\put(172,-2){$t_3=(0,0,1,0,0)$}

\end{picture}}
\end{picture}
\end{center}

\begin{center}
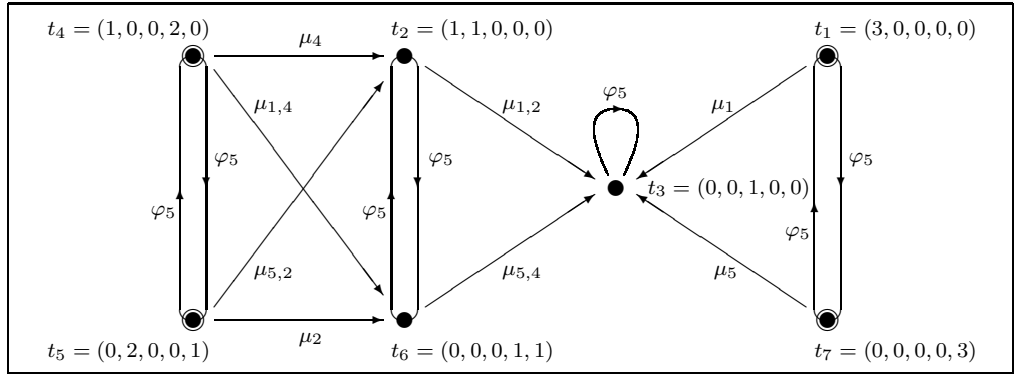
{\bf {Figure 1}} Vertices of $P(G_6,\overline{3})$ with $\mu$-operations and automorphisms.
\linebreak Support vertices are marked by additional circles; \linebreak the bars over the group elements are suppressed.
\end{center}
\vspace{2mm}

Recall that for a group $H$ acting on a set $X,$ the orbit $xH$ of $x\in X$ is the equivalence class of $x$ under $H:$
$$
xH=\{y\in X | y=x\cdot h, \text{for some} ~h\in H\}.
$$
Then the orbits $xH,$ $x\in X,$ form a partition of $X.$

For an integer point $t\in P(G,g_0),$ let us call the non-ordered tuple of its non-zero components $t(g), g\in G_t,$ the multiplicity type of $t.$
For example the point $t=(3,0,1,2,0,1,0^{14})$ of $P(G_{21},\overline{20})$ (not a vertex) is of the multiplicity type $\langle3,2,1,1\rangle.$
As every $\varphi\in Aut(G),$ when applied to some $t\in V(G,g_0),$ rearranges its components without changing their values,
any $s\in tAut(G)$ is of the same multiplicity type as $t$, hence the multiplicity type is the orbit attribute.
However, this does not mean that vertices from different orbits are necessarily of different multiplicity types.

Using Gomory's theorem, Corollary \ref{Cor_Mu-Phi}, Theorem \ref{Th_SupV-Auto}, and Corollary \ref{Cor_Sup-Vert},
we can summarize the description of $V(G,g_0).$

\begin{TT}\label{Th_Vert-Orbits}
The set $V(G,g_0)$ of vertices of any master corner polyhedron $P(G,g_0),$ $g_0\!\neq \!\overline{0},$
is the disjoint union of orbits under the action of $Aut_{g_0}(G).$
All vertices in each orbit are of the same multiplicity type.
For any $t\in V(G,g_0)$ and any $\mu$-operation $\mu$ applicable to $t,$ the map
$t\cdot\varphi\rightarrow\mu_\varphi(t\cdot\varphi),$ $\varphi\in Aut_{g_0}(G),$ sends the orbit $tAut_{g_0}(G)$
vertex-wise onto the orbit $\mu(t)Aut_{g_0}(G).$
Some orbits consist exclusively of support vertices, while the others are free of them.
The set $S(G,g_0)$ of support vertices of $P(G,g_0)$ is the disjoint union of the orbits formed by support vertices.

For the master corner polyhedron $P(G,\overline{0}),$ the analogous assertions hold under $Aut(G)$ acting
on the set $V(G,\overline{0})$ of its vertices.
\end{TT}

As we see in Figure~1, the set of vertices of $P(G_6,\overline{3})$ is partitioned to four orbits under $Aut_{\overline{3}}(G_6)=\{\varepsilon,
\varphi_5\},$ while $S(G_6,\overline{3})$ is the union of two orbits, the very left and the very right ones.

\subsection{Bases}
In this subsection, we describe the bases of the sets of vertices of master corner polyhedra.
These are in some sense the minimal subsets of $V(G,g_0),$ from which one can build all other vertices with the use of the specified tools.
Since extreme rays of all $P(G,g_0)$ are the $g$-axes, $g\in G^+,$ every basis completely determines the polyhedron $P(G,g_0).$

We have at hand two types of instruments that transform vertices to vertices: the automorphisms of the underlying group $G$
that fix the element $g_0$ and the $\mu$-operations. So we can talk about three types of vertex bases:
\begin{itemize}
\item[$B_A$] --- a minimal subset of vertices of a $P(G,g_0)$ such that any other vertex is equivalent to some $b\in B_A;$
\item[$B_S$] --- a minimal subset of vertices of a $P(G,g_0)$ such that any other vertex results by recursive application of some
$\mu$-operations to some $b\in B_S;$
\item[$B_{AS}$] --- a minimal subset of vertices of a $P(G,g_0)$ such that any other vertex results by recursive application of some
automorphisms of $G$ and/or some $\mu$-operations to some $b\in B_{AS}.$
\end{itemize}

It follows from Theorem~\ref{Th_Vert-Orbits} that any system of the orbit representatives for $Aut_{g_0}(G)$ acting on $V(G,g_0)$ is a $B_A$-basis,
hence such basis is unique up to the vertex equivalence.

By the definition of support vertices, the $B_S$-basis is exactly the set of support vertices, thus it is unique.

We have seen in Example 3 that $\{t_1,t_4\}$ is one of the $B_{AS}$-bases for $P(G_6,\overline{3}).$
The $B_{AS}$-bases are described in the next corollary.

\begin{CC}\label{CC_B_AS}
Each $B_{AS}$-basis of the set of vertices of $P(G,g_0)$ is a system of the orbit representatives for $Aut_{g_0}(G)$
acting on the subset of support vertices.
\end{CC}

\noindent{\bf Proof.}
By Gomory's theorem on automorphisms, $V(G,g_0)$ is the union of the orbits under $Aut_{g_0}(G).$
By Theorem~\ref{Th_Vert-Orbits}, each orbit consists of either only support or non-support vertices.
So, each $B_{AS}$-basis must contain exactly one vertex from every orbit consisting of support vertices.
Conversely, each set of such orbit representatives can be taken as a $B_{AS}$-basis. Corollary is proved. $\hfill\square$

It is clear that each $B_A$- and $B_S$-basis contains some $B_{AS}$-subbasis, though the inclusion $B_{AS}\subseteq B_S$ is not necessarily strict.
Introducing of support vertices is evidently advantageous for the polyhedra $P(G,g_0)$ with $|Aut_{g_0}(G)|=1.$
Any $P(G_D,\overline{D-1}),$ $D>3,$ is of this kind since the congruence $k(D-1)\equiv D-1\!\pmod D$ is not fulfilled for any $k,1<k<D.$
Hence, not any two vertices of $P(G_D,\overline{D-1})$ are equivalent, though many of them are not support.
In particular the evidently non-support vertex $s_0$ with the single nonzero component $s_0(g_0)=1$ belongs to every $B_A$-basis
as it is the only vertex of the multiplicity type $\langle1\rangle.$ This yields the strict inclusion $B_{AS}\subset B_A$ for these polyhedra.
The same inclusion holds for all $P(G_D,\overline{0})$ with $D>4,$ since the vertex $t=(1,0^{D-3},1)\in B_A$ equals
$\mu_{\overline{D-2},\overline{1}}(s)$ for the vertex $s=(2,0^{D-4},1,0)$ and thus $t$ is not support.

On the other hand, the strict inclusion $B_{AS}\subset B_S$ holds in particular for all $P(G_D,\overline{r}),$ $r\neq \overline{0},\overline{1},$
with $|Aut_{\overline{r}}(G_D)|>1,$ since for any $\varphi_k\in Aut_{\overline{r}}(G_D),$
$k\neq 1,$ the vertices $(r,0^{D-2})$ and $(0^{k-1},r,0^{D-1-k})$ are both support and equivalent.

Thus, the $B_{AS}$-bases are really the smallest vertex bases for the master corner polyhedra known by now.
Some number characteristics related to vertices of master corner polyhedra $P(G,g_0)$ for all groups $G$ of the order $D<12$
and all right-hand-side elements $g_0$ in the equation (\ref{eq1.1}) are presented in Table \ref{tabular:NumbVertBas} below.
Actually, these are all polyhedra used in \cite{Gom69}. For each $P(G,g_0),$ one can observe there
\begin{itemize}
\item[$\bullet$] $|V(G,g_0)|,$ the number of its vertices,
\item[$\bullet$] $|S(G,g_0)|,$ the number of its support vertices,
\item[$\bullet$] $|B_A|,$ the number of its non-equivalent vertices,
\item[$\bullet$] $|B_{AS}|,$ the cardinality of its $B_{AS}$-bases,
\item[$\bullet$] $|Aut_{g_0}\!(G)|\!-\!1,$ the number of non-identity automorphisms of $G$ that fix $g_0.$
\end{itemize}
The specific values for $G$ and $g_0$ are given in the title of each mini-column. For example
$G_{4,2}$ and $(3,\!0)$ over the last mini-column of the last but one macro-line say that this mini-column refers to the polyhedron $P(G_{4,2},(3,0)),$
where $G_{4,2}$ is the direct sum of the cyclic groups $G_4$ and $G_2$ and $g_0=(\overline{3},\overline{0}).$
The bars over the group elements are suppressed.

{
%\noindent
\begin{table}[htb]
\caption{Vertex characteristics of the $P(G,g_0),$ $|G|<12$} \footnotesize %\scriptsize
\label{tabular:NumbVertBas}
\begin{center}
\begin{tabular}{lcccccccccc}
\toprule
$P(\cdot,\cdot)$    & $G_2,\!0$ & $G_2,\!1$ & $G_3,\!0$ & $G_3,\!2$ & $G_4,\!0$ & $G_4,\!2$ & $G_4,\!3$ & $G_5,\!0$ & $G_5,\!4$ & \\ \hline
                                             %$(D,r)$             & (2,0) & (2,1) & (3,0) & (3,2)
                                             % & (4,0) & (4,2) & (4,3) & (2,2),(0,0) & (2,2),(1,0) \\ \hline
$|V(G,g_0)|$        & 1 & 1 & 3 & 2 & 4 & 3 & 3 & 10 & 5 & \\ %\hline
$|S(G,g_0)|$        & 1 & 1 & 3 & 2 & 4 & 2 & 2 & 8 & 4 & \\ %\hline
$|B_A|$             & 1 & 1 & 2 & 2 & 3 & 2 & 3 & 3 & 5 & \\ %\hline
$|B_{AS}|$          & 1 & 1 & 2 & 2 & 3 & 1 & 2 & 2 & 4 & \\ %\hline
$|Aut_{g_0}\!(G)|\!-\!1$  & 0 & 0 & 1 & 0 & 0 & 1 & 0 & 3 & 0 & \\ %\hline
&  &  &  &  &  &  &  &  & & \\ \midrule [0.07 em]
$P(\cdot,\cdot)$    & $G_6,\!0$ & $G_6,\!3$ & $G_6,\!4$ & $G_6,\!5$ & $G_7,\!0$ & $G_7,\!6$ & $G_8,\!0$ & $G_8,\!4$ & $G_8,\!6$ & $G_8,\!7$ \\
\hline
$|V(G,g_0)|$        & 9 & 7 & 5 & 7 & 23 & 10 & 22 & 9 & 10 & 16 \\ %\hline
$|S(G,g_0)|$        & 7 & 4 & 4 & 4 & 14 & 7 & 14 & 8 & 7 & 9 \\ %\hline
$|B_A|$             & 6 & 4 & 5 & 7 & 5 & 10 & 9 & 4 & 7 & 16 \\ %\hline
$|B_{AS}|$          & 4 & 2 & 4 & 4 & 3 & 7 & 5 & 3 & 4 & 9 \\ %\hline
$|Aut_{g_0}\!(G)|\!-\!1$  & 1 & 1 & 0 & 0 & 5 & 0 & 3 & 3 & 1 & 0 \\ %\hline
&  &  &  &  &  &  &  &  & & \\ \midrule [0.07 em]
$P(\cdot,\cdot)$    & $G_9,\!0$ & $G_9,\!6$ & $G_9,\!8$ & $G_{10},\!0$ & $G_{10},\!5$ & $G_{10},\!8$ & $G_{10},\!9$ & $G_{11},\!0$ & $G_{11},10$ \\ \hline
$|V(G,g_0)|$        & 36 & 14 & 19 & 39 & 29 & 17 & 31 & 85 & 32 \\ %\hline
$|S(G,g_0)|$        & 20 & 10 & 9 & 22 & 16 & 10 & 14 & 40 & 16 \\ %\hline
$|B_A|$             & 7 & 6 & 19 & 11 & 8 & 17 & 31 & 9 & 32 \\ %\hline
$|B_{AS}|$          & 4 & 4 & 9 & 5 & 4 & 10 & 14 & 4 & 16 \\ %\hline
$|Aut_{g_0}\!(G)|\!-\!1$  & 5 & 2 & 0 & 3 & 3 & 0 & 0 & 9 & 0 \\ %\hline
&  &  &  &  &  &  &  &  & \\ \midrule [0.07 em]
$P(\cdot,\cdot)$ & \multicolumn{2}{c}{$G_{2,2},\!(0,\!0)$} & \multicolumn{2}{c}{$G_{2,2},\!(1,\!0)$} & \multicolumn{2}{c}{$G_{4,\!2},\!(0,\!0)$} & \multicolumn{2}{c}{$G_{4,\!2},\!(2,\!0)$} & \multicolumn{2}{c}{$G_{4,\!2},\!(3,\!0)$} \\ \hline
$|V(G,g_0)|$  & \multicolumn{2}{c}{3} & \multicolumn{2}{c}{2} & \multicolumn{2}{c}{9} & \multicolumn{2}{c}{6} & \multicolumn{2}{c}{14} \\ %\hline
$|S(G,g_0)|$  & \multicolumn{2}{c}{3} & \multicolumn{2}{c}{1} & \multicolumn{2}{c}{9} & \multicolumn{2}{c}{6} & \multicolumn{2}{c}{6} \\ %\hline
$|B_A|$ & \multicolumn{2}{c}{3} & \multicolumn{2}{c}{2} & \multicolumn{2}{c}{7} & \multicolumn{2}{c}{4} & \multicolumn{2}{c}{8} \\ %\hline
$|B_{AS}|$ & \multicolumn{2}{c}{3} & \multicolumn{2}{c}{1} & \multicolumn{2}{c}{7} & \multicolumn{2}{c}{4} & \multicolumn{2}{c}{3} \\ %\hline
$|Aut_{g_0}\!(G)|\!-\!1$ & \multicolumn{2}{c}{0} & \multicolumn{2}{c}{0} & \multicolumn{2}{c}{1} & \multicolumn{2}{c}{1} & \multicolumn{2}{c}{1} \\ %\hline
&  &  &  &  &  &  &  &  & \\ \midrule [0.07 em]
$P(\cdot,\cdot)$ & \multicolumn{2}{c}{$G_{2,2,2},\!(0,\!0,\!0)$} & \multicolumn{2}{c}{$G_{2,2,2},\!(1,\!0,\!0)$} & \multicolumn{2}{c}{$G_{3,\!3},\!(0,\!0)$} & \multicolumn{2}{c}{$G_{3,\!3},\!(3,\!0)$} \\ \hline
$|V(G,g_0)|$  & \multicolumn{2}{c}{7} & \multicolumn{2}{c}{8} & \multicolumn{2}{c}{12} & \multicolumn{2}{c}{14} \\ %\hline
$|S(G,g_0)|$  & \multicolumn{2}{c}{7} & \multicolumn{2}{c}{4} & \multicolumn{2}{c}{12} & \multicolumn{2}{c}{3} \\ %\hline
$|B_A|$ & \multicolumn{2}{c}{7} & \multicolumn{2}{c}{8} & \multicolumn{2}{c}{6} & \multicolumn{2}{c}{9} \\ %\hline
$|B_{AS}|$ & \multicolumn{2}{c}{7} & \multicolumn{2}{c}{4} & \multicolumn{2}{c}{6} & \multicolumn{2}{c}{2} \\ %\hline
$|Aut_{g_0}\!(G)|\!-\!1$ & \multicolumn{2}{c}{0} & \multicolumn{2}{c}{0} & \multicolumn{2}{c}{3} & \multicolumn{2}{c}{1} \\ \bottomrule
\end{tabular}
\end{center}
\end{table}
}

More details on these polyhedra, including the lists of vertices and support vertices, can be found in the Appendix.

\section{Adjacency}
By Theorem \ref{Th-MuAdjVert}, $\mu$-operations produce adjacent vertices of master corner polyhedra.
In this section, we show how to construct some complete subgraphs of the graph of the $P(G,g_0)$ applying these operations recursively,
which will help us determine the diameter of this polyhedron.

Let an operation $\mu_{h,f}$ (or $\mu_{h}$) be applied to a vertex $t\in V(G,g_0).$
Let us call $h\in G_t$ the leading element of $t$ and $h+f\in G^+$ (or $t(h)h\in G^+$) the new element of $\mu_{h,f}(t)$ (or of $\mu_{h}(t)$).
Construct a sequence of vertices
\begin{equation}\label{eq4.2}
t=t_0,t_1,t_2,\ldots,t_m=s
\end{equation}
such that $t_i=\mu(t_{i-1}),$ $i=1,\ldots,m,$ for an appropriate $\mu$-operation and the new element of $t_i$ is always the
leading element in constructing $t_{i+1},$ $i=1,\ldots,m-1.$

\begin{TT} \label{Th:CompleteSubgr}
The subgraph of the graph of $P(G,g_0)$ generated by the vertices of the sequence (\ref{eq4.2}) is complete.
\end{TT}

\noindent{\bf Proof.}
By Theorem~\ref{Th-MuAdjVert}, theorem is true for $m=1,$ so we consider $m>1.$
It is sufficient to prove that $s$ is adjacent to $t$ as this would imply that each pair of vertices in (\ref{eq4.2}) are adjacent.
All vertices in (\ref{eq4.2}) are distinct since otherwise there would exist a tuple of integers $u=(u(g),g\in H\subset G_t),$
$u(g)\leq t(g),$ such that $\sum_{g\in H}u(g)g=\overline{0},$ which would contradict Lemma \ref{lem1-Irred} (ii).
Therefore, (\ref{eq4.2}) is a chain in the graph of $P(G,g_0).$

Using the notation from the proof of Theorem~\ref{Th-MuAdjVert}, a collection $\mathcal{F}_{t_1}$ of $D-1$ linearly independent facets of $P(G,g_0)$
that pass through $t_1$ can be chosen in a way to differ from the analogous collection $\mathcal{F}_{t_0}$ 
by containing the trivial facet $t(h_{0})\ge0$ instead of $t(h_{1})\ge0.$
As the leading elements $h_{1},h_{2},\ldots,h_{m-1}$ of $t_1,t_2,\ldots,t_{m-1}$ are their new elements,
the analogous assertions hold true for some collections $\mathcal{F}_{t_2},\ldots,\mathcal{F}_{t_m}.$
They differ from $\mathcal{F}_{t_0}$ by the only one facet: $t(h_{0})\ge0$ instead of
$t(h_{2})\ge0,$ $t(h_{3})\ge0,\ldots,$ $t(h_{m-1}+f_{m-1})\ge0$ with some $f_{m-1}\in G_{t_{m-1}},$ $f_{m-1}\neq h_{m-1},$
or $t(h_{0})\ge0$ instead of $t(t(h_{m-1})h_{m-1})\ge0,$ depending on the type of the $\mu$-operation used at last step.

So, at the end, we obtain that some collection $\mathcal{F}_s$ differs from $\mathcal{F}_t$ by the only one facet,
therefore $s$ is adjacent to $t.$ Theorem is proved. $\hfill\square$

\begin{CC}
The vertex $s_0$ with the components $s_0(g_0)=1,$ $s_0(g)=0$ for all $g\in G^+,$ $g\neq g_0,$ is adjacent to every other vertex of $P(G,g_0).$
\end{CC}

\noindent{\bf Proof.}
At first, notice that $s_0$ is a vertex of $P(G,g_0).$ For arbitrary $t\in V(G,g_0),$ $t\neq s_{0},$ we can build a chain (\ref{eq4.2})
with the end vertex $s=s_0$ in the following way. If $t(h_0)>1$ for some $h_0\in G_t$ then let $t_1=\mu_{h_0}(t),$ so that $h_1=t_{h_0}h_0.$
Otherwise, if $t(g)=1$ for all $g\in G_t,$ set $t_1=t$ and take arbitrary $g\in G_t$ as $h_1.$
In any case, by Lemma \ref{lem1-Irred} (i), $t_1(h_1)=1.$ From this step ahead, we can build the sequence (\ref{eq4.2}) as far as possible,
applying only the operations $\mu_{h,f}$ and taking care only of choosing the leading $h$ elements but not of the $f$ elements.
No matter how we deal forth, we will always obtain a chain with the end vertex $s_0$ since finally all elements of $G_t$ will be added together,
giving $g_0$ as the total sum. Now, it remains to apply Theorem \ref{Th:CompleteSubgr} to complete the proof. $\hfill\square$

Recall that the diameter of a polyhedron is the diameter of its graph.
\begin{TT}
All master corner polyhedra are of diameter 2.
\end{TT}

\noindent{\bf Proof} follows straightforwardly from Corollary 4 as each pair of vertices of $P(G,g_0)$ is connected by a path
of the length $2$ going through $s_0.$ $\hfill\square$

%\newpage
\section{Conclusion}

This paper goes some distance to a better understanding of the vertex-facial structure of the master corner polyhedra.
Any advance in this direction may eventually find implementation in the algorithms to solve the ILP problems.
Our results on the vertex adjacency and on the dependence of nontrivial facets on the vertices they pass through appear to answer this purpose most directly.
The first might be of help to solve the group optimization problem, while the second --- to build facets of corner polyhedra
and then the cutting planes.

However, since corner polyhedra are tightly related to the polyhedra of the ILP problems,
a comprehension of the intriguing vertex structure of the master corner polyhedron might occur of greater importance.
In particular, the very fact that there exists a basic subset of support vertices that often is much smaller in cardinality
compared to the total amount of vertices and even to the number of non-equivalent vertices (see Table 1 and the Appendix) looks promising.
Introducing $\mu$-operations helped us to disclose the non-symmetric relationships between vertices
that complement the symmetric relationships traditionally studied by automorphisms.
If the automorphisms partition the set of vertices to orbits, the new operations connect some pairs of orbits.
The concordance of $\mu$-operations with the subadditivity relations, that play a crucial role in the facet characterization,
and the first computational results move us to believe that hardly any other operation can reduce the bases ever more.

This paper was significantly influenced by the author's recent study of another polyhedron --- the polytope of integer partitions,
which can be regarded as a polyhedron on a set with one partial operation, cf. \cite{Shh05eujc}, \cite{Shh08danvert}, \cite{Shh09danoper}.
Both polyhedra turn out to be rather close relatives.
In particular the nontrivial facets of the partition polytope satisfy a subadditive characterization
similar to (\ref{eq4:ab}) and it also has support vertices.
Numerical data manifest in considerable decrease in the cardinality when going from partitions to vertices
and then to support vertices. One can expect a similar picture for the master corner polyhedron.
What are the numbers of vertices, support vertices, orbits consisting of support vertices ($=|B_{AS}|$),
and the multiplicity types of support vertices --- these are some questions for the future study.

Of special importance is the search for a combinatorial or any other more or less effective characterization of vertices and support vertices.
By now, such characterization is known only for the points inexpressible as a convex combination of two others.
Convex combination is not an easy operation.

\section{Acknowledgements}

It was a surprise to discover that vertices of corner polyhedra had fallen out of the scope
of research in the ILP theory and I am highly grateful to Ralph Gomory, Ellis Johnson, and Jean-Philippe Richard
for the information on the state of the art.

\section{Appendix}

The Appendix presents the tables of vertices of master corner polyhedra $P(G,g_0)$
for all groups $G$ of the order $D<12$ and all right-hand-side elements $g_0$ in the equation~(\ref{eq1.1}).
These tables extend description of their vertex-facet structure presented in \cite{Gom69},
where R.~Gomory listed all their facets, vertices, and the vertex-to-facet incidence matrices.

Each vertex is written as a $(|G|-1)$-dimensional point with zero components suppressed.
The support vertices are marked by the symbol '+' as well as those forming one of the many possible $B_{AS}$-bases.
In the column 'Non-equivalent', we mark by '+' the vertices forming a $B_A$-bases, i.e. one of the many possible sets
of the orbit representatives of $V(G,g_0)$ under the action of $Aut_{g_0}(G).$

\vspace{12mm}
%%%%%%%%%%%%%%%%%%%%%%%%%%%%%%%%%%%%%%%%%%%%%%          $P(G_3,\overline{0})$             %%%%%%%%%%%%%%%%%%%%%%%%%%%%%%%%%%
\begin{table}[htb]
\caption{Vertices of $P(G_3,\overline{0})$}
\footnotesize
\label{tabular:G3-0}
\centering
\begin{tabular}{cc|ccc}
%\hline
\toprule
1 & 2 & Support & Non-equivalent & $B_{AS}$\\ \midrule [0.07 em]
3 &   & + & + & + \\ \hline
1 & 1 & + & + & +  \\ \hline
  & 3 & + &  \\ \bottomrule
\normalsize
\end{tabular}
\end{table}
%%%%%%%%%%%%%%%%%%%%%%%%%%%%%%%%%%%%%%%%%%%%%%          $P(G_3,\overline{0})$             %%%%%%%%%%%%%%%%%%%%%%%%%%%%%%%%%%

%%%%%%%%%%%%%%%%%%%%%%%%%%%%%%%%%%%%%%%%%%%%%%          $P(G_3,\overline{2})$             %%%%%%%%%%%%%%%%%%%%%%%%%%%%%%%%%%
%\newpage
\begin{table}[htb]
\caption{Vertices of $P(G_3,\overline{2})$}
\footnotesize
\label{tabular:G3-2}
\centering
\begin{tabular}{cc|ccc}
\toprule
1 & 2 & Support & Non-equivalent & $B_{AS}$\\ \midrule [0.07 em]
2 &   & + & + & + \\ \hline
  & 1 & + & + & + \\ \bottomrule
\end{tabular}
\end{table}
%%%%%%%%%%%%%%%%%%%%%%%%%%%%%%%%%%         $P(G_3,\overline{2})$

%%%%%%%%%%%%%%%%%%%%%%%%%%%%%%%%%%         $P(G_4,\overline{0})$             %%%%%%%%%%%%%%%%%%%%%%%%%%%%%%%%%%
%\newpage
\begin{table}[htb]
\caption{Vertices of $P(G_4,\overline{0})$} \footnotesize
\label{tabular:G4-0} \centering
\begin{tabular}{ccc|ccc}
\toprule
1 & 2 & 3 & Support & Non-equivalent & $B_{AS}$\\ \midrule [0.07 em]
4 &   &   & + & + & +\\ \hline
  & 2 &   & + & + & +\\ \hline
1 &   & 1 & + & + & +\\ \hline
  &   & 4 & + &  &  \\ \bottomrule
\end{tabular}
\end{table}
%%%%%%%%%%%%%%%%%%%%%%%%%%%%%%%%%%         $P(G_4,\overline{0})$             %%%%%%%%%%%%%%%%%%%%%%%%%%%%%%%%%%

%%%%%%%%%%%%%%%%%%%%%%%%%%%%%%%%%%         $P(G_4,\overline{2})$             %%%%%%%%%%%%%%%%%%%%%%%%%%%%%%%%%%
%\newpage
\begin{table}[htb]
\caption{Vertices of $P(G_4,\overline{2})$}
\footnotesize
\label{tabular:G4-2}
\centering
%\begin{tabular}{|c|c|c|c|c|c|c|c|c|c|}
\begin{tabular}{ccc|ccc}
\toprule
1 & 2 & 3 & Support & Non-equivalent & $B_{AS}$\\ \midrule [0.07 em]
2 &   &   & + & + & + \\ \hline
  & 1 &   &   & + & \\ \hline
  &   & 2 & + & + & + \\ \bottomrule
\end{tabular}
\end{table}
%%%%%%%%%%%%%%%%%%%%%%%%%%%%%%%%%%         $P(G_4,\overline{2})$           %%%%%%%%%%%%%%%%%%%%%%%%%%%%%%%%%%

%%%%%%%%%%%%%%%%%%%%%%%%%%%%%%%%%%         $P(G_4,\overline{3})$           %%%%%%%%%%%%%%%%%%%%%%%%%%%%%%%%%%
%\newpage
\begin{table}[htb]
\caption{Vertices of $P(G_4,\overline{3})$}
\footnotesize
\label{tabular:G4-3}
\centering
\begin{tabular}{ccc|ccc}
\toprule
1 & 2 & 3 & Support & Non-equivalent & $B_{AS}$\\ \midrule [0.07 em]
3 &   &   & + & + & +\\ \hline
1 & 1 &   & + & + & +\\ \hline
  &   & 1 &   & + &  \\ \bottomrule
\end{tabular}
\end{table}
%%%%%%%%%%%%%%%%%%%%%%%%%%%%%%%%%%         $P(G_4,\overline{3})$           %%%%%%%%%%%%%%%%%%%%%%%%%%%%%%%%%%

%%%%%%%%%%%%%%%%%%%%%%%%%%         $P(G_{2,2},(0,0))$           %%%%%%%%%%%%%%%%%%%%%%%%%%%%%%%%%%
\newpage
\begin{table}[h]
\caption{Vertices of $P(G_{2,2},(0,0))$}
\footnotesize
\label{tabular:G2,2-0,0} \centering
\begin{tabular}{ccc|ccc}
\toprule
(1,0) & (0,1) & (1,1) & Support & Non-equivalent & $B_{AS}$\\ \midrule [0.07 em]
2 &   &   & + & + & + \\ \hline
  & 2 &   & + & + & + \\ \hline
  &   & 2 & + & + & + \\ \bottomrule
\end{tabular}
\end{table}
%%%%%%%%%%%%%%%%%%%%%%%%%%         $P(G_{2,2},(0,0))$           %%%%%%%%%%%%%%%%%%%%%%%%%%%%%%%%%%

%%%%%%%%%%%%%%%%%%%%%%%%%%%%%%%%%%%%%%%%%%%%%%%%%%%%%%%%         $P(G_{2,2},(1,0))$           %%%%%%%%%%%%%%%%%%%%%%%%%%%%%%%%%%
\newpage
\begin{table}[h]
\caption{Vertices of $P(G_{2,2},(1,0))$}
\footnotesize
\label{tabular:G2,2-1,0} \centering
\begin{tabular}{ccc|ccc}
\toprule
(1,0) & (0,1) & (1,1) & Support & Non-equivalent & $B_{AS}$\\ \midrule [0.07 em]
1 &   &   &   & + &   \\ \hline
  & 1 & 1 & + & + & + \\ \bottomrule
\end{tabular}
\end{table}
%%%%%%%%%%%%%%%%%%%%%%%%%%%%%%%%%%%%%%%%%%%%%%%%%%%%%%%         $P(G_{2,2},(1,0))$           %%%%%%%%%%%%%%%%%%%%%%%%%%%%%%%%%%

%%%%%%%%%%%%%%%%%%%%%%%%%%%%%%%%%%         $P(G_5,\overline{0})$          %%%%%%%%%%%%%%%%%%%%%%%%%%%%%%%%%%
%\newpage
\begin{table}[h]
\caption{Vertices of $P(G_5,\overline{0})$}
\footnotesize
\label{tabular:G5-0} \centering
\begin{tabular}{cccc|ccc}
\toprule
1 & 2 & 3 & 4 & Support & Non-equivalent & $B_{AS}$\\ \midrule [0.07 em]
5 &   &   &   &  +  & + & +  \\ \hline
1 & 2 &   &   &  +  & + & +  \\ \hline
  & 5 &   &   &  +  &   & \\ \hline
2 &   & 1 &   &  +  &   & \\ \hline
  & 1 & 1 &   &     & + & \\ \hline
  &   & 5 &   &  +  &   & \\ \hline
1 &   &   & 1 &  +  &   & \\ \hline
  &   & 1 & 1 &  +  &   & \\ \hline
  & 1 &   & 2 &  +  &   & \\ \hline
  &   &   & 5 &  +  &   & \\ \bottomrule
\end{tabular}
\end{table}
%%%%%%%%%%%%%%%%%%%%%%%%%%%%%%%%%%         $P(G_5,\overline{0})$          %%%%%%%%%%%%%%%%%%%%%%%%%%%%%%%%%%

%%%%%%%%%%%%%%%%%%%%%%%%%%%%%%%%%%         $P(G_5,\overline{4})$          %%%%%%%%%%%%%%%%%%%%%%%%%%%%%%%%%%
%\newpage
\begin{table}[h]
\caption{Vertices of $P(G_5,\overline{4})$}
\footnotesize
\label{tabular:G5-4} \centering
\begin{tabular}{cccc|ccc}
\toprule
1 & 2 & 3 & 4 & Support & Non-equivalent & $B_{AS}$\\ \midrule [0.07 em]
4 &   &   &   &  +  & + & +  \\ \hline
  & 2 &   &   &  +  & + & +  \\ \hline
1 &   & 1 &   &  +  & + & +  \\ \hline
  &   & 3 &   &  +  & + & +  \\ \hline
  &   &   & 1 &     & + &    \\ \bottomrule
\end{tabular}
\end{table}
%%%%%%%%%%%%%%%%%%%%%%%%%%%%%%%%%%         $P(G_5,\overline{4})$          %%%%%%%%%%%%%%%%%%%%%%%%%%%%%%%%%%

%%%%%%%%%%%%%%%%%%%%%%%%%%%%%%%%%%         $P(G_6,\overline{0})$          %%%%%%%%%%%%%%%%%%%%%%%%%%%%%%%%%%
\newpage
\begin{table}[h]
\caption{Vertices of $P(G_6,\overline{0})$}
\footnotesize
\label{tabular:G6-0} \centering
\begin{tabular}{ccccc|ccc}
\toprule
1 & 2 & 3 & 4 & 5 & Support & Non-equivalent & $B_{AS}$\\ \midrule [0.07 em]
6 &   &   &   &   &  +  &  +   & + \\ \hline
2 &   &   & 1 &   &  +  &  +   & + \\ \hline
1 &   &   &   & 1 &     &   +  &\\ \hline
  & 3 &   &   &   &  +  &   +  & + \\ \hline
  & 1 &   & 1 &   &     &   +  &\\ \hline
  & 1 &   &   & 2 &  +  &    &\\ \hline
  &   & 2 &   &   &  +  &   +  & + \\ \hline
  &   &   & 3 &   &  +  &    &\\ \hline
  &   &   &   & 6 &  +  &    & \\ \bottomrule
\end{tabular}
\end{table}
%%%%%%%%%%%%%%%%%%%%%%%%%%%%%%%%%%         $P(G_6,\overline{0})$          %%%%%%%%%%%%%%%%%%%%%%%%%%%%%%%%%%

%%%%%%%%%%%%%%%%%%%%%%%%%%%%%%%%%%         $P(G_6,\overline{3})$          %%%%%%%%%%%%%%%%%%%%%%%%%%%%%%%%%%
%\newpage
\begin{table}[h]
\caption{Vertices of $P(G_6,\overline{3})$}
\footnotesize
\label{tabular:G6-3} \centering %\begin{tabular}{|c|c|c|c|c|c|c|c|c|c|}
\begin{tabular}{ccccc|ccc}
\toprule
1 & 2 & 3 & 4 & 5 & Support & Non-equivalent & $B_{AS}$\\ \midrule [0.07 em]
3 &   &   &   &   & + & +  & + \\ \hline
1 & 1 &   &   &   &   &  + &   \\ \hline
  &   & 1 &   &   &   &  + &   \\ \hline
1 &   &   & 2 &   & + &  &   \\ \hline
  & 2 &   &   & 1 & + & +  & + \\ \hline
  &   &   & 1 & 1 &   &  &   \\ \hline
  &   &   &   & 3 & + &  &   \\ \bottomrule
\end{tabular}
\end{table}
%%%%%%%%%%%%%%%%%%%%%%%%%%%%%%%%%%         $P(G_6,\overline{3})$          %%%%%%%%%%%%%%%%%%%%%%%%%%%%%%%%%%

%%%%%%%%%%%%%%%%%%%%%%%%%%%%%%%%%%         $P(G_6,\overline{4})$          %%%%%%%%%%%%%%%%%%%%%%%%%%%%%%%%%%
%\newpage
\begin{table}[h]
\caption{Vertices of $P(G_6,\overline{4})$}
\footnotesize
\label{tabular:G6-4}
\centering                                                             %\begin{tabular}{|c|c|c|c|c|c|c|c|c|c|}
\begin{tabular}{ccccc|ccc}
\toprule
$t_1$ & 2 & 3 & 4 & 5 & Support & Non-equivalent & $B_{AS}$\\ \midrule [0.07 em]
4 &   &   &   &   &  +  & + & + \\ \hline
  & 2 &   &   &   &  +  & + & + \\ \hline
1 &   & 1 &   &   &  +  & + & + \\ \hline
  &   &   & 1 &   &     & + &   \\ \hline
  &   &   &   & 2 &  +  & + & + \\ \bottomrule
\end{tabular}
\end{table}
%%%%%%%%%%%%%%%%%%%%%%%%%%%%%%%%%%         $P(G_6,\overline{4})$          %%%%%%%%%%%%%%%%%%%%%%%%%%%%%%%%%%

%%%%%%%%%%%%%%%%%%%%%%%%%%%%%%%%%%         $P(G_6,\overline{5})$          %%%%%%%%%%%%%%%%%%%%%%%%%%%%%%%%%%
\newpage
\begin{table}[h]
\caption{Vertices of $P(G_6,\overline{5})$}
\footnotesize
\label{tabular:G6-5}
\centering
%\begin{tabular}{|c|c|c|c|c|c|c|c|c|c|}
\begin{tabular}{ccccc|ccc}
\toprule
$t_1$ & 2 & 3 & 4 & 5 & Support & Non-equivalent & $B_{AS}$\\ \midrule [0.07 em]
5 &   &   &   &   & + &  +  & + \\ \hline
2 &   & 1 &   &   & + &  +  & + \\ \hline
1 & 2 &   &   &   & + &   + & + \\ \hline
1 &   &   & 1 &   &   &   + &   \\ \hline
  & 1 & 1 &   &   &   &   + &   \\ \hline
  &   & 1 & 2 &   & + &   + & + \\ \hline
  &   &   &   & 1 &   &   + &   \\ \bottomrule
\end{tabular}
\end{table}
%%%%%%%%%%%%%%%%%%%%%%%%%%%%%%%%%%         $P(G_6,\overline{5})$          %%%%%%%%%%%%%%%%%%%%%%%%%%%%%%%%%%

%%%%%%%%%%%%%%%%%%%%%%%%%%%%%%%%%%         $P(G_7,\overline{0})$          %%%%%%%%%%%%%%%%%%%%%%%%%%%%%%%%%%
%\newpage
\begin{table}[h]
\caption{Vertices of $P(G_7,\overline{0})$}
\footnotesize
\label{tabular:G7-0}
\centering
%\begin{tabular}{|c|c|c|c|c|c|c|c|c|c|}
\begin{tabular}{cccccc|ccccc}
\toprule
1 & 2 & 3 & 4 & 5 & 6 & Support & Non-equivalent & $B_{AS}$\\ \midrule [0.07 em]
7 &   &   &   &   &   & + &  +   & + \\ \hline
3 &   &   & 1 &   &   & + &  +   & + \\ \hline
2 &   &   &   & 1 &   &   &   +  &   \\ \hline
1 & 3 &   &   &   &   & + &    &   \\ \hline
1 & 1 &   & 1 &   &   & + & +    & + \\ \hline
1 &   & 2 &   &   &   &   &   &   \\ \hline
1 &   &   &   &   & 1 &   &  +   &   \\ \hline
  & 7 &   &   &   &   & + &    &   \\ \hline
  & 2 & 1 &   &   &   &   &    &   \\ \hline
  & 1 &   & 3 &   &   & + &    &   \\ \hline
  & 1 &   &   & 1 &   &   &    &   \\ \hline
  & 1 &   &   &   & 2 & + &    &   \\ \hline
  &   & 7 &   &   &   & + &    &   \\ \hline
  &   & 3 &   & 1 &   & + &    &   \\ \hline
  &   & 1 & 1 &   &   &   &    &   \\ \hline
  &   & 1 &   & 1 & 1 & + &    &   \\ \hline
  &   & 1 &   &   & 3 & + &    &   \\ \hline
  &   &   & 7 &   &   & + &    &   \\ \hline
  &   &   & 2 &   & 1 &   &   &   \\ \hline
  &   &   & 1 & 2 &   &   &   &   \\ \hline
  &   &   &   & 7 &   & + &   &   \\ \hline
  &   &   &   & 3 & 1 & + &   &   \\ \hline
  &   &   &   &   & 7 & + &   &   \\ \bottomrule
\end{tabular}
\end{table}
%%%%%%%%%%%%%%%%%%%%%%%%%%%%%%%%%%         $P(G_7,\overline{0})$          %%%%%%%%%%%%%%%%%%%%%%%%%%%%%%%%%%

%%%%%%%%%%%%%%%%%%%%%%%%%%%%%%%%%%         $P(G_7,\overline{6})$          %%%%%%%%%%%%%%%%%%%%%%%%%%%%%%%%%%
\newpage
\begin{table}[h]
\caption{Vertices of $P(G_7,\overline{6})$}
\footnotesize
\label{tabular:G7-6}
\centering
\begin{tabular}{cccccc|ccc}
\toprule
1 & 2 & 3 & 4 & 5 & 6 & Support & Non-equivalent & $B_{AS}$\\ \midrule [0.07 em]
6 &   &   &   &   &   & + & +  & + \\ \hline
2 &   &   & 1 &   &   & + & +  & + \\ \hline
1 &   &   &   & 1 &   &   & +  \\ \hline
  & 3 &   &   &   &   & + &  + & + \\ \hline
  & 1 &   & 1 &   &   &   & +  \\ \hline
  &   & 2 &   &   &   & + & +  & + \\ \hline
  &   &   & 5 &   &   & + &  + & + \\ \hline
  &   &   & 2 & 1 &   & + &  + & + \\ \hline
  &   &   &   & 4 &   & + &  + & + \\ \hline
  &   &   &   &   & 1 &   &  + \\ \bottomrule
\end{tabular}
\end{table}
%%%%%%%%%%%%%%%%%%%%%%%%%%%%%%%%%%         $P(G_7,\overline{6})$          %%%%%%%%%%%%%%%%%%%%%%%%%%%%%%%%%%

%%%%%%%%%%%%%%%%%%%%%%%%%%%%%%%%%%         $P(G_8,\overline{0})$          %%%%%%%%%%%%%%%%%%%%%%%%%%%%%%%%%%
%\newpage
\begin{table}[h]
\caption{Vertices of $P(G_8,\overline{0})$}
\footnotesize
\label{tabular:G8-0}
\centering
\begin{tabular}{ccccccc|ccc}
\toprule
1 & 2 & 3 & 4 & 5 & 6 & 7 & Support & Non-equivalent & $B_{AS}$\\ \midrule [0.07 em]
8 & 0 & 0 & 0 & 0 & 0 & 0 &  +  &  +   & +  \\ \hline
3 & 0 & 0 & 0 & 1 & 0 & 0 &  +  &   +  & +  \\ \hline
2 & 0 & 2 & 0 & 0 & 0 & 0 &  +  &   +  & +  \\ \hline
2 & 0 & 0 & 0 & 0 & 1 & 0 &     &   +  & \\ \hline
1 & 1 & 0 & 0 & 1 & 0 & 0 &  +  &   +  & +  \\ \hline
1 & 0 & 0 & 0 & 3 & 0 & 0 &  +  &    & \\ \hline
1 & 0 & 0 & 0 & 0 & 0 & 1 &     &  +   & \\ \hline
0 & 4 & 0 & 0 & 0 & 0 & 0 &  +  &   +  & +  \\ \hline
0 & 1 & 2 & 0 & 0 & 0 & 0 &     &    & \\ \hline
0 & 1 & 0 & 0 & 0 & 1 & 0 &     &  +   & \\ \hline
0 & 1 & 0 & 0 & 0 & 0 & 2 &     &   & \\ \hline
0 & 0 & 8 & 0 & 0 & 0 & 0 &  +  &   & \\ \hline
0 & 0 & 3 & 0 & 0 & 0 & 1 &  +  &   & \\ \hline
0 & 0 & 1 & 0 & 1 & 0 & 0 &     &   & \\ \hline
0 & 0 & 1 & 0 & 0 & 1 & 1 &  +  &   & \\ \hline
0 & 0 & 1 & 0 & 0 & 0 & 3 &  +  &   & \\ \hline
0 & 0 & 0 & 2 & 0 & 0 & 0 &     &  +   & \\ \hline
0 & 0 & 0 & 0 & 8 & 0 & 0 &  +  &   & \\ \hline
0 & 0 & 0 & 0 & 2 & 1 & 0 &     &   & \\ \hline
0 & 0 & 0 & 0 & 2 & 0 & 2 &  +  &   & \\ \hline
0 & 0 & 0 & 0 & 0 & 4 & 0 &  +  &   & \\ \hline
0 & 0 & 0 & 0 & 0 & 0 & 8 &  +  &   & \\ \bottomrule
\end{tabular}
\end{table}
%%%%%%%%%%%%%%%%%%%%%%%%%%%%%%%%%%         $P(G_8,\overline{0})$          %%%%%%%%%%%%%%%%%%%%%%%%%%%%%%%%%%

%%%%%%%%%%%%%%%%%%%%%%%%%%%%%%%%%%         $P(G_8,\overline{4})$          %%%%%%%%%%%%%%%%%%%%%%%%%%%%%%%%%%
\newpage
\begin{table}[h]
\caption{Vertices of $P(G_8,\overline{4})$}
\footnotesize
\label{tabular:G8-4}
\centering
\begin{tabular}{ccccccc|ccc}
\toprule
1 & 2 & 3 & 4 & 5 & 6 & 7 & Support & Non-equivalent & $B_{AS}$\\ \midrule [0.07 em]
4 &   &   &   &   &   &   & + &  + & + \\ \hline
1 &   & 1 &   &   &   &   & + & +  & + \\ \hline
  & 2 &   &   &   &   &   & + & +  & + \\ \hline
  &   & 4 &   &   &   &   & + &  &   \\ \hline
  &   &   & 1 &   &   &   &   & +  &   \\ \hline
  &   &   &   & 4 &   &   & + & &   \\ \hline
  &   &   &   & 1 &   & 1 & + & &   \\ \hline
  &   &   &   &   & 2 &   & + & &   \\ \hline
  &   &   &   &   &   & 4 & + & &   \\ \bottomrule
\end{tabular}
\end{table}
%%%%%%%%%%%%%%%%%%%%%%%%%%%%%%%%%%         $P(G_8,\overline{4})$          %%%%%%%%%%%%%%%%%%%%%%%%%%%%%%%%%%

%%%%%%%%%%%%%%%%%%%%%%%%%%%%%%%%%%         $P(G_8,\overline{6})$          %%%%%%%%%%%%%%%%%%%%%%%%%%%%%%%%%%
%\newpage
\begin{table}[h]
\caption{Vertices of $P(G_8,\overline{6})$}
\footnotesize
\label{tabular:G8-6}
\centering
\begin{tabular}{ccccccc|ccc}
\toprule
1 & 2 & 3 & 4 & 5 & 6 & 7 & Support & Non-equivalent & $B_{AS}$\\ \midrule [0.07 em]
6 &   &   &   &   &   &   & + & +  & + \\ \hline
2 &   &   & 1 &   &   &   & + & +  & + \\ \hline
1 &   &   &   & 1 &   &   &   & +  &   \\ \hline
  & 3 &   &   &   &   &   & + & +  & + \\ \hline
  & 1 &   & 1 &   &   &   &   & +  &   \\ \hline
  &   & 2 &   &   &   &   & + & +  & + \\ \hline
  &   &   & 1 & 2 &   &   & + &  &   \\ \hline
  &   &   &   & 6 &   &   & + &  &   \\ \hline
  &   &   &   &   & 1 &   &   & +  &   \\ \hline
  &   &   &   &   &   & 2 & + &  &   \\ \bottomrule
\end{tabular}
\end{table}
%%%%%%%%%%%%%%%%%%%%%%%%%%%%%%%%%%         $P(G_8,\overline{6})$          %%%%%%%%%%%%%%%%%%%%%%%%%%%%%%%%%%

%%%%%%%%%%%%%%%%%%%%%%%%%%%%%%%%%%         $P(G_8,\overline{7})$          %%%%%%%%%%%%%%%%%%%%%%%%%%%%%%%%%%
%\vspace{-8mm}
\newpage
\begin{table}[h]
\caption{Vertices of $P(G_8,\overline{7})$}
\footnotesize
\label{tabular:G8-7}
\centering
\begin{tabular}{ccccccc|ccc}
\toprule
1 & 2 & 3 & 4 & 5 & 6 & 7 & Support & Non-equivalent & $B_{AS}$\\ \midrule [0.07 em]
7 &   &   &   &   &   &   & + & + & + \\ \hline
3 &   &   & 1 &   &   &   & + & + & + \\ \hline
2 &   &   &   & 1 &   &   &   & + &  \\ \hline
1 & 3 &   &   &   &   &   & + & + & + \\ \hline
1 & 1 &   & 1 &   &   &   & + & + & + \\ \hline
1 &   & 2 &   &   &   &   & + & + & + \\ \hline
1 &   &   &   &   & 1 &   &   & + &  \\ \hline
  & 2 & 1 &   &   &   &   &   & + &  \\ \hline
  & 1 &   &   & 1 &   &   &   & + &  \\ \hline
  &   & 5 &   &   &   &   & + & + & + \\ \hline
  &   & 1 & 1 &   &   &   &   & + &  \\ \hline
  &   & 1 &   &   & 2 &   &   & + &  \\ \hline
  &   &   & 1 & 1 & 1 &   & + & + & + \\ \hline
  &   &   &   & 3 &   &   & + & + & + \\ \hline
  &   &   &   & 1 & 3 &   & + & + & + \\ \hline
  &   &   &   &   &   & 1 &   & + &  \\ \bottomrule
\end{tabular}
\end{table}
%%%%%%%%%%%%%%%%%%%%%%%%%%%%%%%%%%         $P(G_8,\overline{7})$          %%%%%%%%%%%%%%%%%%%%%%%%%%%%%%%%%%

%%%%%%%%%%%%%%%%%%%%%%%%%%%%%%%%%%         $P(G_{4,2},(0,0))$           %%%%%%%%%%%%%%%%%%%%%%%%%%%%%%%%%%
%\newpage
\begin{table}[h] \caption{Vertices of $P(G_{4,2},(0,0))$} \footnotesize
\label{tabular:G4,2-0,0} \centering
\begin{tabular}{ccccccc|ccc}
\toprule
(1,0) & (2,0) & (3,0) & (0,1) & (1,1) & (2,1) & (3,1) & Support & Non-equivalent & $B_{AS}$\\ \midrule [0.07 em]
4 & 0 & 0 & 0 & 0 & 0 & 0 &  +  & + & + \\ \hline
1 & 0 & 1 & 0 & 0 & 0 & 0 &  +  & + & + \\ \hline
0 & 2 & 0 & 0 & 0 & 0 & 0 &  +  & + & + \\ \hline
0 & 0 & 4 & 0 & 0 & 0 & 0 &  +  & + & + \\ \hline
0 & 0 & 0 & 2 & 0 & 0 & 0 &  +  & + & + \\ \hline
0 & 0 & 0 & 0 & 4 & 0 & 0 &  +  &    & \\ \hline
0 & 0 & 0 & 0 & 0 & 2 & 0 &  +  & + & + \\ \hline
0 & 0 & 0 & 0 & 1 & 0 & 1 &  +  & + & + \\ \hline
0 & 0 & 0 & 0 & 0 & 0 & 4 &  +  &    & \\ \bottomrule
\end{tabular}
\end{table}
%%%%%%%%%%%%%%%%%%%%%%%%%%%%%%%%%%         $P(G_{4,2},(0,0))$           %%%%%%%%%%%%%%%%%%%%%%%%%%%%%%%%%%

%%%%%%%%%%%%%%%%%%%%%%%%%%%%%%%%%%         $P(G_{4,2},(2,0))$           %%%%%%%%%%%%%%%%%%%%%%%%%%%%%%%%%%
%\newpage
\begin{table}[h]
\caption{Vertices of $P(G_{4,2},(2,0))$} \footnotesize
\label{tabular:G4,2-2,0} \centering
\begin{tabular}{ccccccc|ccc}
\toprule
(1,0) & (2,0) & (3,0) & (0,1) & (1,1) & (2,1) & (3,1) & Support & Non-equivalent & $B_{AS}$\\ \midrule [0.07 em]
2 & 0 & 0 & 0 & 0 & 0 & 0 &  +  &  +  & + \\ \hline
0 & 1 & 0 & 0 & 0 & 0 & 0 &  +  &   + & + \\ \hline
0 & 0 & 2 & 0 & 0 & 0 & 0 &  +  &   & \\ \hline
0 & 0 & 0 & 1 & 0 & 1 & 0 &  +  &  +  & + \\ \hline
0 & 0 & 0 & 0 & 2 & 0 & 0 &  +  &  +  & + \\ \hline
0 & 0 & 0 & 0 & 0 & 0 & 2 &  +  &   & \\ \bottomrule
\end{tabular}
\end{table}
%%%%%%%%%%%%%%%%%%%%%%%%%%%%%%%%%%         $P(G_{4,2},(2,0))$           %%%%%%%%%%%%%%%%%%%%%%%%%%%%%%%%%%

%%%%%%%%%%%%%%%%%%%%%%%%%%%%%%%%%%         $P(G_{4,2},(3,0))$           %%%%%%%%%%%%%%%%%%%%%%%%%%%%%%%%%%
\newpage
\begin{table}[h]
\caption{Vertices of $P(G_{4,2},(3,0))$} \footnotesize
\label{tabular:G4,2-3,0} \centering
\begin{tabular}{ccccccc|ccc} \toprule
(1,0) & (2,0) & (3,0) & (0,1) & (1,1) & (2,1) & (3,1) & Support & Non-equivalent & $B_{AS}$\\ \midrule [0.07 em]
3 & 0 & 0 & 0 & 0 & 0 & 0 &  +  & + & + \\ \hline
1 & 1 & 0 & 0 & 0 & 0 & 0 &     & + & \\ \hline
1 & 0 & 0 & 1 & 0 & 1 & 0 &  +  & + & + \\ \hline
1 & 0 & 0 & 0 & 2 & 0 & 0 &     & + & \\ \hline
1 & 0 & 0 & 0 & 0 & 0 & 2 &  +  & + & + \\ \hline
0 & 1 & 0 & 1 & 1 & 0 & 0 &  +  & + & + \\ \hline
0 & 1 & 0 & 0 & 0 & 1 & 1 &  +  & + & + \\ \hline
0 & 0 & 1 & 0 & 0 & 0 & 0 &     & + & \\ \hline
0 & 0 & 0 & 1 & 3 & 0 & 0 &  +  & + & + \\ \hline
0 & 0 & 0 & 1 & 0 & 0 & 1 &     & + & \\ \hline
0 & 0 & 0 & 0 & 1 & 1 & 0 &     & + & \\ \hline
0 & 0 & 0 & 0 & 0 & 1 & 3 &  +  & + & + \\ \bottomrule
\end{tabular}
\end{table}
%%%%%%%%%%%%%%%%%%%%%%%%%%%%%%%%%%         $P(G_{4,2},(3,0))$           %%%%%%%%%%%%%%%%%%%%%%%%%%%%%%%%%%

%%%%%%%%%%%%%%%%%%%%%%%%%%%%%%%%%%         $P(G_{4,2},(0,1))$           %%%%%%%%%%%%%%%%%%%%%%%%%%%%%%%%%%
%\newpage
\begin{table}[h]
\caption{Vertices of $P(G_{4,2},(0,1))$} \footnotesize
\label{tabular:G4,2-0,1} \centering
\begin{tabular}{ccccccc|ccc} \toprule
(1,0) & (2,0) & (3,0) & (0,1) & (1,1) & (2,1) & (3,1) & Support & Non-equivalent & $B_{AS}$\\ \midrule [0.07 em]
3 & 0 & 0 & 0 & 1 & 0 & 0 &  +  & + & + \\ \hline
2 & 0 & 0 & 0 & 0 & 1 & 0 &     & + & \\ \hline
1 & 1 & 0 & 0 & 1 & 0 & 0 &  +  & + & + \\ \hline
1 & 0 & 0 & 0 & 3 & 0 & 0 &  +  & + & + \\ \hline
1 & 0 & 0 & 0 & 0 & 0 & 1 &     & + & \\ \hline
0 & 1 & 1 & 0 & 0 & 0 & 1 &  +  &    & \\ \hline
0 & 1 & 0 & 0 & 0 & 1 & 0 &     & + & \\ \hline
0 & 0 & 3 & 0 & 0 & 0 & 1 &  +  &   & \\ \hline
0 & 0 & 2 & 0 & 0 & 1 & 0 &     &   & \\ \hline
0 & 0 & 1 & 0 & 1 & 0 & 0 &     &   & \\ \hline
0 & 0 & 1 & 0 & 0 & 0 & 3 &  +  &   & \\ \hline
0 & 0 & 0 & 1 & 0 & 0 & 0 &     & + & \\ \hline
0 & 0 & 0 & 0 & 2 & 1 & 0 &     & + & \\ \hline
0 & 0 & 0 & 0 & 0 & 1 & 2 &     &    &  \\ \bottomrule
\end{tabular}
\end{table}
%%%%%%%%%%%%%%%%%%%%%%%%%%%%%%%%%%         $P(G_{4,2},(0,1))$           %%%%%%%%%%%%%%%%%%%%%%%%%%%%%%%%%%

%%%%%%%%%%%%%%%%%%%%%%%%%%%%%%%%%%         $P(G_{2,2,2},(0,0,0))$           %%%%%%%%%%%%%%%%%%%%%%%%%%%%%%%%%%
\newpage
\begin{table}[h]
\caption{Vertices of $P(G_{2,2,2},(0,0,0))$} \footnotesize
\label{tabular:G2,2,2-0,0,0} \centering
\begin{tabular}{ccccccc|ccc} \toprule
(1,0,0) & (0,1,0) & (1,1,0) & (0,0,1) & (1,0,1) & (0,1,1) & (1,1,1) & Support & Non-equivalent & $B_{AS}$\\ \midrule [0.07 em]
2 & 0 & 0 & 0 & 0 & 0 & 0 &  +  & + & + \\ \hline
0 & 2 & 0 & 0 & 0 & 0 & 0 &  +  & + & + \\ \hline
0 & 0 & 2 & 0 & 0 & 0 & 0 &  +  & + & + \\ \hline
0 & 0 & 0 & 2 & 0 & 0 & 0 &  +  & + & + \\ \hline
0 & 0 & 0 & 0 & 2 & 0 & 0 &  +  & + & + \\ \hline
0 & 0 & 0 & 0 & 0 & 2 & 0 &  +  & + & + \\ \hline
0 & 0 & 0 & 0 & 0 & 0 & 2 &  +  & + & + \\ \bottomrule
\end{tabular}
\end{table}
%%%%%%%%%%%%%%%%%%%%%%%%%%%%%%%%%%         $P(G_{2,2,2},(0,0,0))$           %%%%%%%%%%%%%%%%%%%%%%%%%%%%%%%%%%

%%%%%%%%%%%%%%%%%%%%%%%%%%%%%%%%%%         $P(G_{2,2,2},(1,0,0))$           %%%%%%%%%%%%%%%%%%%%%%%%%%%%%%%%%%
%\newpage
\begin{table}[h]
\caption{Vertices of $P(G_{2,2,2},(1,0,0))$} \footnotesize
\label{tabular:G2,2,2-1,0,0} \centering
\begin{tabular}{ccccccc|ccc} \toprule
(1,0,0) & (0,1,0) & (1,1,0) & (0,0,1) & (1,0,1) & (0,1,1) & (1,1,1) & Support & Non-equivalent & $B_{AS}$\\ \midrule [0.07 em]
1 & 0 & 0 & 0 & 0 & 0 & 0 &     & + & \\ \hline
0 & 1 & 1 & 0 & 0 & 0 & 0 &     & + & \\ \hline
0 & 1 & 0 & 1 & 0 & 0 & 1 &  +  & + & + \\ \hline
0 & 1 & 0 & 0 & 1 & 1 & 0 &  +  & + & + \\ \hline
0 & 0 & 1 & 1 & 0 & 1 & 0 &  +  & + & + \\ \hline
0 & 0 & 1 & 0 & 1 & 0 & 1 &  +  & + & + \\ \hline
0 & 0 & 0 & 1 & 1 & 0 & 0 &     & + & \\ \hline
0 & 0 & 0 & 0 & 0 & 1 & 1 &     & + & \\ \bottomrule
\end{tabular}
\end{table}
%%%%%%%%%%%%%%%%%%%%%%%%%%%%%%%%%%         $P(G_{2,2,2},(1,0,0))$           %%%%%%%%%%%%%%%%%%%%%%%%%%%%%%%%%%

%%%%%%%%%%%%%%%%%%%%%%%%%%%%%%%%%%         $P(G_9,\overline{0})$          %%%%%%%%%%%%%%%%%%%%%%%%%%%%%%%%%%
\newpage
\begin{table}[h]
\caption{Vertices of $P(G_9,\overline{0})$}
\footnotesize
\label{tabular:G9-0}
\centering
\begin{tabular}{cccccccc|ccc}
\toprule
1 & 2 & 3 & 4 & 5 & 6 & 7 & 8 & Support & Non-equivalent & $B_{AS}$\\ \midrule [0.07 em]
9 & 0 & 0 & 0 & 0 & 0 & 0 & 0 &  +  & + & + \\ \hline
4 & 0 & 0 & 0 & 1 & 0 & 0 & 0 &  +  & + & + \\ \hline
3 & 0 & 0 & 0 & 0 & 1 & 0 & 0 &     & + & \\ \hline
2 & 0 & 0 & 0 & 0 & 0 & 1 & 0 &     & + & \\ \hline
1 & 4 & 0 & 0 & 0 & 0 & 0 & 0 &  +  &    & \\ \hline
1 & 1 & 0 & 0 & 0 & 1 & 0 & 0 &  +  & + & + \\ \hline
1 & 0 & 0 & 2 & 0 & 0 & 0 & 0 &     &   & \\ \hline
1 & 0 & 1 & 0 & 1 & 0 & 0 & 0 &  +  &   & \\ \hline
1 & 0 & 0 & 0 & 0 & 0 & 0 & 1 &     &   & \\ \hline
0 & 9 & 0 & 0 & 0 & 0 & 0 & 0 &  +  &   & \\ \hline
0 & 3 & 1 & 0 & 0 & 0 & 0 & 0 &     &   & \\ \hline
0 & 2 & 0 & 0 & 1 & 0 & 0 & 0 &     &   & \\ \hline
0 & 1 & 1 & 1 & 0 & 0 & 0 & 0 &  +  &   & \\ \hline
0 & 1 & 0 & 4 & 0 & 0 & 0 & 0 &  +  &   & \\ \hline
0 & 1 & 0 & 0 & 0 & 0 & 1 & 0 &     &   & \\ \hline
0 & 1 & 0 & 0 & 0 & 0 & 0 & 2 &     &   & \\ \hline
0 & 0 & 3 & 0 & 0 & 0 & 0 & 0 &  +  &  +  & + \\ \hline
0 & 0 & 1 & 0 & 3 & 0 & 0 & 0 &     &    & \\ \hline
0 & 0 & 1 & 0 & 0 & 1 & 0 & 0 &     &  +   & \\ \hline
0 & 0 & 1 & 0 & 0 & 0 & 1 & 1 &  +  &   & \\ \hline
0 & 0 & 1 & 0 & 0 & 0 & 0 & 3 &     &  & \\ \hline
0 & 0 & 0 & 9 & 0 & 0 & 0 & 0 &  +  &   & \\ \hline
0 & 0 & 0 & 3 & 0 & 1 & 0 & 0 &     &  & \\ \hline
0 & 0 & 0 & 1 & 1 & 0 & 0 & 0 &     &  & \\ \hline
0 & 0 & 0 & 1 & 0 & 1 & 0 & 1 &  +  &   & \\ \hline
0 & 0 & 0 & 1 & 0 & 0 & 2 & 0 &     &  & \\ \hline
0 & 0 & 0 & 1 & 0 & 0 & 0 & 4 &  +  &   & \\ \hline
0 & 0 & 0 & 0 & 9 & 0 & 0 & 0 &  +  &   & \\ \hline
0 & 0 & 0 & 0 & 4 & 0 & 1 & 0 &  +  &   & \\ \hline
0 & 0 & 0 & 0 & 2 & 0 & 0 & 1 &     &  & \\ \hline
0 & 0 & 0 & 0 & 1 & 1 & 1 & 0 &  +  &   & \\ \hline
0 & 0 & 0 & 0 & 0 & 3 & 0 & 0 &  +  &   & \\ \hline
0 & 0 & 0 & 0 & 0 & 1 & 3 & 0 &     &  & \\ \hline
0 & 0 & 0 & 0 & 0 & 0 & 9 & 0 &  +  &   & \\ \hline
0 & 0 & 0 & 0 & 0 & 0 & 4 & 1 &  +  &    & \\ \hline
0 & 0 & 0 & 0 & 0 & 0 & 0 & 9 &  +  &    & \\ \bottomrule
\end{tabular}
\end{table}
%%%%%%%%%%%%%%%%%%%%%%%%%%%%%%%%%%         $P(G_9,\overline{0})$          %%%%%%%%%%%%%%%%%%%%%%%%%%%%%%%%%%

%%%%%%%%%%%%%%%%%%%%%%%%%%%%%%%%%%         $P(G_9,\overline{6})$          %%%%%%%%%%%%%%%%%%%%%%%%%%%%%%%%%%
\newpage
\begin{table}[h]
\caption{Vertices of $P(G_9,\overline{6})$}
\footnotesize
\label{tabular:G9-6}
\centering
\begin{tabular}{cccccccc|ccc}
\toprule
1 & 2 & 3 & 4 & 5 & 6 & 7 & 8 & Support & Non-equivalent & $B_{AS}$\\ \midrule [0.07 em]
6 & 0 & 0 & 0 & 0 & 0 & 0 & 0 &  +  & + & + \\ \hline
2 & 0 & 0 & 1 & 0 & 0 & 0 & 0 &  +  & + & + \\ \hline
1 & 0 & 0 & 0 & 1 & 0 & 0 & 0 &     & + & \\ \hline
1 & 0 & 0 & 0 & 0 & 0 & 2 & 0 &  +  &   & \\ \hline
0 & 3 & 0 & 0 & 0 & 0 & 0 & 0 &  +  & + & + \\ \hline
0 & 1 & 0 & 1 & 0 & 0 & 0 & 0 &     &    & \\ \hline
0 & 0 & 2 & 0 & 0 & 0 & 0 & 0 &  +  & + & + \\ \hline
0 & 0 & 0 & 6 & 0 & 0 & 0 & 0 &  +  &    & \\ \hline
0 & 0 & 0 & 2 & 0 & 0 & 1 & 0 &  +  &    & \\ \hline
0 & 0 & 0 & 0 & 3 & 0 & 0 & 0 &  +  &    & \\ \hline
0 & 0 & 0 & 0 & 0 & 1 & 0 & 0 &     & + & \\ \hline
0 & 0 & 0 & 0 & 0 & 0 & 6 & 0 &  +  &    & \\ \hline
0 & 0 & 0 & 0 & 0 & 0 & 1 & 1 &     &    & \\ \hline
0 & 0 & 0 & 0 & 0 & 0 & 0 & 3 &  +  &    & \\ \bottomrule
\end{tabular}
\end{table}
%%%%%%%%%%%%%%%%%%%%%%%%%%%%%%%%%%         $P(G_9,\overline{6})$          %%%%%%%%%%%%%%%%%%%%%%%%%%%%%%%%%%

%%%%%%%%%%%%%%%%%%%%%%%%%%%%%%%%%%         $P(G_9,\overline{8})$          %%%%%%%%%%%%%%%%%%%%%%%%%%%%%%%%%%
%\newpage
\begin{table}[h]
\caption{Vertices of $P(G_9,\overline{8})$}
\footnotesize
\label{tabular:G9-8}
\centering
\begin{tabular}{cccccccc|ccc}
\toprule
1 & 2 & 3 & 4 & 5 & 6 & 7 & 8 & Support & Non-equivalent & $B_{AS}$\\ \midrule [0.07 em]
8 & 0 & 0 & 0 & 0 & 0 & 0 & 0 &  +  & + & + \\ \hline
3 & 0 & 0 & 0 & 1 & 0 & 0 & 0 &  +  & + & + \\ \hline
2 & 0 & 2 & 0 & 0 & 0 & 0 & 0 &  +  & + & + \\ \hline
2 & 0 & 0 & 0 & 0 & 1 & 0 & 0 &     & + & \\ \hline
1 & 1 & 0 & 0 & 1 & 0 & 0 & 0 &  +  & + & + \\ \hline
1 & 0 & 0 & 0 & 0 & 0 & 1 & 0 &     & + & \\ \hline
0 & 4 & 0 & 0 & 0 & 0 & 0 & 0 &  +  & + & + \\ \hline
0 & 1 & 2 & 0 & 0 & 0 & 0 & 0 &     & + & \\ \hline
0 & 1 & 0 & 0 & 3 & 0 & 0 & 0 &  +  & + & + \\ \hline
0 & 1 & 0 & 0 & 0 & 1 & 0 & 0 &     & + & \\ \hline
0 & 0 & 1 & 0 & 1 & 0 & 0 & 0 &     & + & \\ \hline
0 & 0 & 1 & 0 & 0 & 0 & 2 & 0 &     & + & \\ \hline
0 & 0 & 0 & 2 & 0 & 0 & 0 & 0 &     & + & \\ \hline
0 & 0 & 0 & 0 & 7 & 0 & 0 & 0 &  +  & + & + \\ \hline
0 & 0 & 0 & 0 & 2 & 0 & 1 & 0 &     & + & \\ \hline
0 & 0 & 0 & 0 & 1 & 2 & 0 & 0 &     & + & \\ \hline
0 & 0 & 0 & 0 & 0 & 2 & 2 & 0 &  +  & + & + \\ \hline
0 & 0 & 0 & 0 & 0 & 0 & 5 & 0 &  +  & + & + \\ \hline
0 & 0 & 0 & 0 & 0 & 0 & 0 & 1 &     & + & \\ \bottomrule
\end{tabular}
\end{table}
%%%%%%%%%%%%%%%%%%%%%%%%%%%%%%%%%%         $P(G_9,\overline{8})$          %%%%%%%%%%%%%%%%%%%%%%%%%%%%%%%%%%

%%%%%%%%%%%%%%%%%%%%%%%%%%%%%%%%%%         $P(G_{3,3},(0,0))$           %%%%%%%%%%%%%%%%%%%%%%%%%%%%%%%%%%
\newpage
\begin{table}[h]
\caption{Vertices of $P(G_{3,3},(0,0))$} \footnotesize
\label{tabular:G3,3-0,0} \centering
\begin{tabular}{cccccccc|ccc} \toprule
(1,0) & (2,0) & (0,1) & (1,1) & (2,1) & (0,2) & (1,2) & (2,2) & Support & Non-equivalent & $B_{AS}$\\ \midrule [0.07 em]
3 & 0 & 0 & 0 & 0 & 0 & 0 & 0 &  +  & + & + \\ \hline
1 & 1 & 0 & 0 & 0 & 0 & 0 & 0 &  +  & + & + \\ \hline
0 & 3 & 0 & 0 & 0 & 0 & 0 & 0 &  +  &    & \\ \hline
0 & 0 & 3 & 0 & 0 & 0 & 0 & 0 &  +  & + & + \\ \hline
0 & 0 & 1 & 0 & 0 & 1 & 0 & 0 &  +  & + & + \\ \hline
0 & 0 & 0 & 3 & 0 & 0 & 0 & 0 &  +  & + & + \\ \hline
0 & 0 & 0 & 0 & 3 & 0 & 0 & 0 &  +  &    & \\ \hline
0 & 0 & 0 & 1 & 0 & 0 & 0 & 1 &  +  & + & + \\ \hline
0 & 0 & 0 & 0 & 1 & 0 & 1 & 0 &  +  &    & \\ \hline
0 & 0 & 0 & 0 & 0 & 3 & 0 & 0 &  +  &    & \\ \hline
0 & 0 & 0 & 0 & 0 & 0 & 3 & 0 &  +  &    & \\ \hline
0 & 0 & 0 & 0 & 0 & 0 & 0 & 3 &  +  &    & \\ \bottomrule
\end{tabular}
\end{table}
%%%%%%%%%%%%%%%%%%%%%%%%%%%%%%%%%%         $P(G_{3,3},(0,0))$           %%%%%%%%%%%%%%%%%%%%%%%%%%%%%%%%%%

%%%%%%%%%%%%%%%%%%%%%%%%%%%%%%%%%%         $P(G_{3,3},(1,0))$           %%%%%%%%%%%%%%%%%%%%%%%%%%%%%%%%%%
%\newpage
\begin{table}[h]
\caption{Vertices of $P(G_{3,3},(1,0))$} \footnotesize
\label{tabular:G3,3-1,0} \centering
\begin{tabular}{cccccccc|ccc} \toprule
(1,0) & (2,0) & (0,1) & (1,1) & (2,1) & (0,2) & (1,2) & (2,2) & Support & Non-equivalent & $B_{AS}$\\ \midrule [0.07 em]
1 & 0 & 0 & 0 & 0 & 0 & 0 & 0 &     & + & \\ \hline
0 & 2 & 0 & 0 & 0 & 0 & 0 & 0 &     & + & \\ \hline
0 & 0 & 2 & 1 & 0 & 0 & 0 & 0 &     & + & \\ \hline
0 & 0 & 2 & 0 & 0 & 0 & 0 & 2 &  +  & + & + \\ \hline
0 & 0 & 1 & 0 & 2 & 0 & 0 & 0 &     & + & \\ \hline
0 & 0 & 1 & 0 & 0 & 0 & 1 & 0 &     & + & \\ \hline
0 & 0 & 0 & 2 & 1 & 0 & 0 & 0 &     & + & \\ \hline
0 & 0 & 0 & 2 & 0 & 0 & 2 & 0 &  +  & + & + \\ \hline
0 & 0 & 0 & 1 & 0 & 1 & 0 & 0 &     &    & \\ \hline
0 & 0 & 0 & 0 & 2 & 2 & 0 & 0 &  +  &    & \\ \hline
0 & 0 & 0 & 0 & 1 & 0 & 0 & 1 &     & + & \\ \hline
0 & 0 & 0 & 0 & 0 & 2 & 1 & 0 &     &    & \\ \hline
0 & 0 & 0 & 0 & 0 & 1 & 0 & 2 &     &    & \\ \hline
0 & 0 & 0 & 0 & 0 & 0 & 2 & 1 &     &    & \\ \bottomrule
\end{tabular}
\end{table}
%%%%%%%%%%%%%%%%%%%%%%%%%%%%%%%%%%         $P(G_{3,3},(1,0))$           %%%%%%%%%%%%%%%%%%%%%%%%%%%%%%%%%%

%%%%%%%%%%%%%%%%%%%%%%%%%%%%%%%%%%         $P(G_10,\overline{0})$          %%%%%%%%%%%%%%%%%%%%%%%%%%%%%%%%%%
\newpage
\begin{table}[h]
\caption{Vertices of $P(G_{10},\overline{0})$}
\footnotesize
\label{tabular:G10-0}
\centering
\begin{tabular}{ccccccccc|ccc}
\toprule
1 & 2 & 3 & 4 & 5 & 6 & 7 & 8 & 9 & Support & Non-equivalent & $B_{AS}$\\ \midrule [0.07 em]
10 & 0 & 0 & 0 & 0 & 0 & 0 & 0 & 0 &  +  & + & + \\ \hline
4 & 0 & 0 & 0 & 0 & 1 & 0 & 0 & 0 &  +  & + & + \\ \hline
3 & 0 & 0 & 0 & 0 & 0 & 1 & 0 & 0 &     & + & \\ \hline
2 & 0 & 0 & 2 & 0 & 0 & 0 & 0 & 0 &  +  & + & + \\ \hline
2 & 0 & 0 & 0 & 0 & 0 & 0 & 1 & 0 &     & + & \\ \hline
1 & 0 & 3 & 0 & 0 & 0 & 0 & 0 & 0 &     &    & \\ \hline
1 & 1 & 0 & 0 & 0 & 0 & 1 & 0 & 0 &  +  & + & + \\ \hline
1 & 0 & 1 & 0 & 0 & 1 & 0 & 0 & 0 &  +  &    & \\ \hline
1 & 0 & 0 & 0 & 0 & 0 & 0 & 0 & 1 &     & + & \\ \hline
0 & 5 & 0 & 0 & 0 & 0 & 0 & 0 & 0 &  +  & + & + \\ \hline
0 & 2 & 2 & 0 & 0 & 0 & 0 & 0 & 0 &  +  &    & \\ \hline
0 & 2 & 0 & 0 & 0 & 1 & 0 & 0 & 0 &     & + & \\ \hline
0 & 1 & 0 & 2 & 0 & 0 & 0 & 0 & 0 &     &    & \\ \hline
0 & 1 & 0 & 0 & 0 & 0 & 4 & 0 & 0 &  +  &    & \\ \hline
0 & 1 & 0 & 0 & 0 & 0 & 0 & 1 & 0 &     & + & \\ \hline
0 & 1 & 0 & 0 & 0 & 0 & 0 & 0 & 2 &     &    & \\ \hline
0 & 0 & 10 & 0 & 0 & 0 & 0 & 0 & 0 &  +  &    & \\ \hline
0 & 0 & 4 & 0 & 0 & 0 & 0 & 1 & 0 &  +  &   & \\ \hline
0 & 0 & 2 & 1 & 0 & 0 & 0 & 0 & 0 &     &  & \\ \hline
0 & 0 & 1 & 0 & 0 & 0 & 1 & 0 & 0 &     &  & \\ \hline
0 & 0 & 1 & 0 & 0 & 0 & 0 & 1 & 1 &  +  &   & \\ \hline
0 & 0 & 1 & 0 & 0 & 0 & 0 & 0 & 3 &     &  & \\ \hline
0 & 0 & 0 & 5 & 0 & 0 & 0 & 0 & 0 &  +  &   & \\ \hline
0 & 0 & 0 & 1 & 0 & 1 & 0 & 0 & 0 &     &  & \\ \hline
0 & 0 & 0 & 1 & 0 & 0 & 1 & 0 & 1 &  +  &   & \\ \hline
0 & 0 & 0 & 1 & 0 & 0 & 0 & 2 & 0 &     &  & \\ \hline
0 & 0 & 0 & 1 & 0 & 0 & 0 & 0 & 4 &  +  &   & \\ \hline
0 & 0 & 0 & 0 & 2 & 0 & 0 & 0 & 0 &     & + & \\ \hline
0 & 0 & 0 & 0 & 1 & 1 & 0 & 0 & 1 &  +  &   & \\ \hline
0 & 0 & 0 & 0 & 1 & 0 & 1 & 1 & 0 &  +  &   & \\ \hline
0 & 0 & 0 & 0 & 0 & 5 & 0 & 0 & 0 &  +  &   & \\ \hline
0 & 0 & 0 & 0 & 0 & 2 & 0 & 1 & 0 &     &   & \\ \hline
0 & 0 & 0 & 0 & 0 & 2 & 0 & 0 & 2 &  +  &   & \\ \hline
0 & 0 & 0 & 0 & 0 & 1 & 2 & 0 & 0 &     &   & \\ \hline
0 & 0 & 0 & 0 & 0 & 0 & 10 & 0 & 0 &  +  &    & \\ \hline
0 & 0 & 0 & 0 & 0 & 0 & 3 & 0 & 1 &     &  & \\ \hline
0 & 0 & 0 & 0 & 0 & 0 & 2 & 2 & 0 &  +  &   & \\ \hline
0 & 0 & 0 & 0 & 0 & 0 & 0 & 5 & 0 &  +  &   & \\ \hline
0 & 0 & 0 & 0 & 0 & 0 & 0 & 0 & 10 &  +  &    & \\ \bottomrule
\end{tabular}
\end{table}
%%%%%%%%%%%%%%%%%%%%%%%%%%%%%%%%%%         $P(G_10,\overline{0})$          %%%%%%%%%%%%%%%%%%%%%%%%%%%%%%%%%%

%%%%%%%%%%%%%%%%%%%%%%%%%%%%%%%%%%         $P(G_10,\overline{5})$          %%%%%%%%%%%%%%%%%%%%%%%%%%%%%%%%%%
\newpage
\begin{table}[h]
\caption{Vertices of $P(G_{10},\overline{5})$} \footnotesize
\label{tabular:G10-5} \centering
\begin{tabular}{ccccccccc|ccc}
\toprule
1 & 2 & 3 & 4 & 5 & 6 & 7 & 8 & 9 & Support & Non-equivalent & $B_{AS}$\\ \midrule [0.07 em]
5 &   &   &   &   &   &   &   &   & + &  +  & + \\ \hline
2 &   & 1 &   &   &   &   &   &   & + &  +  & + \\ \hline
1 & 2 &   &   &   &   &   &   &   &   &  +  \\ \hline
1 &   &   & 1 &   &   &   &   &   &   &  +  \\ \hline
1 &   &   &   &   & 4 &   &   &   & + &   \\ \hline
1 &   &   &   &   & 1 &   & 1 &   & + &  \\ \hline
1 &   &   &   &   &   & 2 &   &   & + &   \\ \hline
1 &   &   &   &   &   &   & 3 &   &   &  \\ \hline
  & 4 &   &   &   &   & 1 &   &   & + & +   & + \\ \hline
  & 3 &   &   &   &   &   &   & 1 &   &  +  \\ \hline
  & 1 & 1 &   &   &   &   &   &   &   &  \\ \hline
  & 1 &   & 1 &   &   &   &   & 1 & + & \\ \hline
  & 1 &   &   &   & 1 & 1 &   &   & + & \\ \hline
  &   & 5 &   &   &   &   &   &   & + &  \\ \hline
  &   & 2 &   &   &   &   &   & 1 & + &  \\ \hline
  &   & 1 & 3 &   &   &   &   &   &   &  \\ \hline
  &   & 1 & 1 &   &   &   & 1 &   & + & +   & + \\ \hline
  &   & 1 &   &   &   &   & 4 &   & + &   \\ \hline
  &   & 1 &   &   & 2 &   &   &   &   &   \\ \hline
  &   &   & 4 &   &   &   &   & 1 & + &   \\ \hline
  &   &   & 2 &   &   & 1 &   &   &   &   \\ \hline
  &   &   &   & 1 &   &   &   &   &   &  +  \\ \hline
  &   &   &   &   & 3 & 1 &   &   &   &  \\ \hline
  &   &   &   &   & 1 &   &   & 1 &   &  \\ \hline
  &   &   &   &   &   & 5 &   &   & + &   \\ \hline
  &   &   &   &   &   & 1 & 1 &   &   &  \\ \hline
  &   &   &   &   &   & 1 &   & 2 & + &  \\ \hline
  &   &   &   &   &   &   & 2 & 1 &   &  \\ \hline
  &   &   &   &   &   &   &   & 5 & + &   \\ \bottomrule
\end{tabular}
\centering
\end{table}
%%%%%%%%%%%%%%%%%%%%%%%%%%%%%%%%%%         $P(G_10,\overline{5})$          %%%%%%%%%%%%%%%%%%%%%%%%%%%%%%%%%%

%%%%%%%%%%%%%%%%%%%%%%%%%%%%%%%%%%         $P(G_10,\overline{8})$          %%%%%%%%%%%%%%%%%%%%%%%%%%%%%%%%%%
\newpage
\begin{table}[h]
\caption{Vertices of $P(G_{10},\overline{8})$}
\footnotesize
\label{tabular:G10-8}
\centering
\begin{tabular}{ccccccccc|ccc}
\toprule
1 & 2 & 3 & 4 & 5 & 6 & 7 & 8 & 9 & Support & Non-equivalent & $B_{AS}$\\ \midrule [0.07 em]
8 & 0 & 0 & 0 & 0 & 0 & 0 & 0 & 0 &  +  & + & + \\ \hline
3 & 0 & 0 & 0 & 1 & 0 & 0 & 0 & 0 &  +  & + & + \\ \hline
2 & 0 & 2 & 0 & 0 & 0 & 0 & 0 & 0 &  +  & + & + \\ \hline
2 & 0 & 0 & 0 & 0 & 1 & 0 & 0 & 0 &     & + & \\ \hline
1 & 1 & 0 & 0 & 1 & 0 & 0 & 0 & 0 &  +  & + & + \\ \hline
1 & 0 & 0 & 0 & 0 & 0 & 1 & 0 & 0 &     & + & \\ \hline
0 & 4 & 0 & 0 & 0 & 0 & 0 & 0 & 0 &  +  & + & + \\ \hline
0 & 1 & 2 & 0 & 0 & 0 & 0 & 0 & 0 &     & + & \\ \hline
0 & 1 & 0 & 0 & 0 & 1 & 0 & 0 & 0 &     & + & \\ \hline
0 & 0 & 6 & 0 & 0 & 0 & 0 & 0 & 0 &  +  & + & + \\ \hline
0 & 0 & 1 & 0 & 1 & 0 & 0 & 0 & 0 &     & + & \\ \hline
0 & 0 & 0 & 2 & 0 & 0 & 0 & 0 & 0 &     & + & \\ \hline
0 & 0 & 0 & 0 & 1 & 1 & 1 & 0 & 0 &  +  & + & + \\ \hline
0 & 0 & 0 & 0 & 0 & 3 & 0 & 0 & 0 &  +  & + & + \\ \hline
0 & 0 & 0 & 0 & 0 & 0 & 4 & 0 & 0 &  +  & + & + \\ \hline
0 & 0 & 0 & 0 & 0 & 0 & 0 & 1 & 0 &     & + & \\ \hline
0 & 0 & 0 & 0 & 0 & 0 & 0 & 0 & 2 &  +  & + & + \\ \bottomrule
\end{tabular}
\end{table}
%%%%%%%%%%%%%%%%%%%%%%%%%%%%%%%%%%         $P(G_10,\overline{8})$          %%%%%%%%%%%%%%%%%%%%%%%%%%%%%%%%%%

%%%%%%%%%%%%%%%%%%%%%%%%%%%%%%%%%%         $P(G_10,\overline{9})$          %%%%%%%%%%%%%%%%%%%%%%%%%%%%%%%%%%
\newpage
\begin{table}[h]
\caption{Vertices of $P(G_{10},\overline{9})$}
\footnotesize
\label{tabular:G10-9}
\centering
\begin{tabular}{ccccccccc|ccc}
\toprule
1 & 2 & 3 & 4 & 5 & 6 & 7 & 8 & 9 & Support & Non-equivalent & $B_{AS}$\\ \midrule [0.07 em]
9 & 0 & 0 & 0 & 0 & 0 & 0 & 0 & 0 &  +  & + & + \\ \hline
4 & 0 & 0 & 0 & 1 & 0 & 0 & 0 & 0 &  +  & + & + \\ \hline
3 & 0 & 0 & 0 & 0 & 1 & 0 & 0 & 0 &     & + & \\ \hline
2 & 0 & 0 & 0 & 0 & 0 & 1 & 0 & 0 &     & + & \\ \hline
1 & 4 & 0 & 0 & 0 & 0 & 0 & 0 & 0 &  +  & + & + \\ \hline
1 & 1 & 0 & 0 & 0 & 1 & 0 & 0 & 0 &  +  & + & + \\ \hline
1 & 0 & 1 & 0 & 1 & 0 & 0 & 0 & 0 &  +  & + & + \\ \hline
1 & 0 & 0 & 2 & 0 & 0 & 0 & 0 & 0 &     & + & \\ \hline
1 & 0 & 0 & 0 & 0 & 3 & 0 & 0 & 0 &     & + & \\ \hline
1 & 0 & 0 & 0 & 0 & 0 & 4 & 0 & 0 &  +  & + & + \\ \hline
1 & 0 & 0 & 0 & 0 & 0 & 0 & 1 & 0 &     & + & \\ \hline
0 & 3 & 1 & 0 & 0 & 0 & 0 & 0 & 0 &     & + & \\ \hline
0 & 2 & 0 & 0 & 1 & 0 & 0 & 0 & 0 &     & + & \\ \hline
0 & 1 & 1 & 1 & 0 & 0 & 0 & 0 & 0 &  +  & + & + \\ \hline
0 & 1 & 0 & 0 & 0 & 0 & 1 & 0 & 0 &     & + & \\ \hline
0 & 0 & 3 & 0 & 0 & 0 & 0 & 0 & 0 &  +  & + & + \\ \hline
0 & 0 & 1 & 4 & 0 & 0 & 0 & 0 & 0 &  +  & + & + \\ \hline
0 & 0 & 1 & 0 & 0 & 1 & 0 & 0 & 0 &     & + & \\ \hline
0 & 0 & 1 & 0 & 0 & 0 & 0 & 2 & 0 &     & + & \\ \hline
0 & 0 & 0 & 3 & 0 & 0 & 1 & 0 & 0 &     & + & \\ \hline
0 & 0 & 0 & 1 & 1 & 0 & 0 & 0 & 0 &     & + & \\ \hline
0 & 0 & 0 & 1 & 0 & 0 & 1 & 1 & 0 &  +  & + & + \\ \hline
0 & 0 & 0 & 0 & 1 & 4 & 0 & 0 & 0 &  +  & + & + \\ \hline
0 & 0 & 0 & 0 & 1 & 1 & 0 & 1 & 0 &  +  & + & + \\ \hline
0 & 0 & 0 & 0 & 1 & 0 & 2 & 0 & 0 &     & + & \\ \hline
0 & 0 & 0 & 0 & 1 & 0 & 0 & 3 & 0 &     & + & \\ \hline
0 & 0 & 0 & 0 & 0 & 2 & 1 & 0 & 0 &     & + & \\ \hline
0 & 0 & 0 & 0 & 0 & 0 & 7 & 0 & 0 &  +  & + & + \\ \hline
0 & 0 & 0 & 0 & 0 & 0 & 3 & 1 & 0 &     & + & \\ \hline
0 & 0 & 0 & 0 & 0 & 0 & 1 & 4 & 0 &  +  & + & + \\ \hline
0 & 0 & 0 & 0 & 0 & 0 & 0 & 0 & 1 &     & + & \\ \bottomrule
\end{tabular}
\end{table}
%%%%%%%%%%%%%%%%%%%%%%%%%%%%%%%%%%         $P(G_10,\overline{9})$          %%%%%%%%%%%%%%%%%%%%%%%%%%%%%%%%%%

%%%%%%%%%%%%%%%%%%%%%%%%%%%%%%%%%%         $P(G_11,\overline{10})$          %%%%%%%%%%%%%%%%%%%%%%%%%%%%%%%%%%
\newpage
\begin{table}[h]
\caption{Vertices of $P(G_{11},\overline{10})$}
\footnotesize
\label{tabular:G10-9}
\centering
\begin{tabular}{cccccccccc|ccc}
\toprule
1 & 2 & 3 & 4 & 5 & 6 & 7 & 8 & 9 & 10 & Support & Non-equivalent & $B_{AS}$\\ \midrule [0.07 em]
10 & 0 & 0 & 0 & 0 & 0 & 0 & 0 & 0 & 0 &  +  & + & + \\ \hline
4 & 0 & 0 & 0 & 0 & 1 & 0 & 0 & 0 & 0 &  +  & + & + \\ \hline
3 & 0 & 0 & 0 & 0 & 0 & 1 & 0 & 0 & 0 &     & + & \\ \hline
2 & 0 & 0 & 0 & 0 & 0 & 0 & 1 & 0 & 0 &     & + & \\ \hline
2 & 0 & 0 & 2 & 0 & 0 & 0 & 0 & 0 & 0 &  +  & + & + \\ \hline
1 & 1 & 0 & 0 & 0 & 0 & 1 & 0 & 0 & 0 &  +  & + & + \\ \hline
1 & 0 & 3 & 0 & 0 & 0 & 0 & 0 & 0 & 0 &  +  & + & + \\ \hline
1 & 0 & 1 & 0 & 0 & 1 & 0 & 0 & 0 & 0 &  +  & + & + \\ \hline
1 & 0 & 0 & 0 & 0 & 0 & 0 & 0 & 1 & 0 &     & + & \\ \hline
0 & 5 & 0 & 0 & 0 & 0 & 0 & 0 & 0 & 0 &  +  & + & + \\ \hline
0 & 2 & 2 & 0 & 0 & 0 & 0 & 0 & 0 & 0 &  +  & + & + \\ \hline
0 & 2 & 0 & 0 & 0 & 1 & 0 & 0 & 0 & 0 &     & + & \\ \hline
0 & 1 & 0 & 2 & 0 & 0 & 0 & 0 & 0 & 0 &     & + & \\ \hline
0 & 1 & 0 & 0 & 0 & 0 & 0 & 1 & 0 & 0 &     & + & \\ \hline
0 & 0 & 7 & 0 & 0 & 0 & 0 & 0 & 0 & 0 &  +  & + & + \\ \hline
0 & 0 & 2 & 1 & 0 & 0 & 0 & 0 & 0 & 0 &     & + & \\ \hline
0 & 0 & 1 & 0 & 0 & 3 & 0 & 0 & 0 & 0 &     & + & \\ \hline
0 & 0 & 1 & 0 & 0 & 0 & 1 & 0 & 0 & 0 &     & + & \\ \hline
0 & 0 & 1 & 0 & 0 & 0 & 0 & 0 & 2 & 0 &  +  & + & + \\ \hline
0 & 0 & 0 & 8 & 0 & 0 & 0 & 0 & 0 & 0 &  +  & + & + \\ \hline
0 & 0 & 0 & 3 & 0 & 0 & 0 & 0 & 1 & 0 &  +  & + & + \\ \hline
0 & 0 & 0 & 1 & 0 & 1 & 0 & 0 & 0 & 0 &     & + & \\ \hline
0 & 0 & 0 & 1 & 0 & 0 & 0 & 1 & 1 & 0 &  +  & + & + \\ \hline
0 & 0 & 0 & 0 & 2 & 0 & 0 & 0 & 0 & 0 &     & + & \\ \hline
0 & 0 & 0 & 0 & 0 & 9 & 0 & 0 & 0 & 0 &  +  & + & + \\ \hline
0 & 0 & 0 & 0 & 0 & 4 & 0 & 1 & 0 & 0 &  +  & + & + \\ \hline
0 & 0 & 0 & 0 & 0 & 2 & 0 & 0 & 1 & 0 &     & + & \\ \hline
0 & 0 & 0 & 0 & 0 & 1 & 1 & 1 & 0 & 0 &  +  & + & + \\ \hline
0 & 0 & 0 & 0 & 0 & 0 & 3 & 0 & 0 & 0 &  +  & + & + \\ \hline
0 & 0 & 0 & 0 & 0 & 0 & 0 & 4 & 0 & 0 &  +  & + & + \\ \hline
0 & 0 & 0 & 0 & 0 & 0 & 0 & 0 & 6 & 0 &  +  & + & + \\ \hline
0 & 0 & 0 & 0 & 0 & 0 & 0 & 0 & 0 & 1 &     & + & \\ \bottomrule
\end{tabular}
\end{table}
%%%%%%%%%%%%%%%%%%%%%%%%%%%%%%%%%%         $P(G_11,\overline{10})$          %%%%%%%%%%%%%%%%%%%%%%%%%%%%%%%%%%

\renewcommand{\baselinestretch}{1}

\end{document}